%% file: adaptive_hss.tex
\documentclass[onefignum,onetabnum]{siamart171218}

\input{tex/preamble.tex}

\begin{document}

  \input{tex/front_matter.tex}

  \input{tex/introduction.tex}

  \input{tex/hss.tex}

\input{tex/parallel_algorithm.tex}

  \input{tex/experiments.tex}

  \input{tex/conclusion.tex}

\input{tex/acknowledgment.tex}

  \input{tex/appendix.tex}

  \bibliographystyle{siamplain}
  \bibliography{adaptive_hss}

\end{document}

%% file: tex/preamble.tex
\usepackage{float}
\usepackage{fancyvrb}

\usepackage{lineno,hyperref}
\usepackage{verbatim}
\modulolinenumbers[5]

\usepackage[algo2e,ruled,boxed,vlined,linesnumbered]{algorithm2e}
\SetKwRepeat{Do}{do}{while}

\usepackage{amsmath,amssymb,amsfonts}

\newcommand{\eps}[0]{\varepsilon}
\newcommand{\R}[0]{\mathbb{R}}
\newcommand{\E}[0]{\mathbb{E}}
\renewcommand{\P}[0]{\mathbb{P}}

\newcommand{\parens}[1]{\left( #1 \right)}
\newcommand{\brackets}[1]{\left[ #1 \right]}
\newcommand{\braces}[1]{\left\{ #1 \right\}}
\newcommand{\abs}[1]{\left| #1 \right|}
\newcommand{\norm}[1]{\left|\left| #1 \right|\right|}

\usepackage{todonotes}
\usepackage{multirow}

\newsiamremark{remark}{Remark}

\usepackage{tikz}
\usetikzlibrary{decorations.pathreplacing}
\tikzstyle{cnode} = [draw, circle,scale=1.0]
\tikzstyle{level 1} = [level distance=.2\textwidth, sibling distance=.5\textwidth]
\tikzstyle{level 2} = [sibling distance=.25\textwidth]
\tikzstyle{level 3} = [sibling distance=.15\textwidth]

\newcommand{\xslcomment}[1]{\textcolor{red}{\bf [#1 --Sherry]}}
\newcommand{\chgcomment}[1]{\textcolor{green}{\bf [#1 --Chris]}}

\newcommand{\ignore}[1]{}
\newcommand{\dbl}{{\textsc{Doubling}}}
\newcommand{\inc}{{\textsc{Incrementing}}}
\newcommand{\PC}{{\texttt{PARTIALLY\_COMPRESSED}}}

%% file: tex/front_matter.tex
\title{Matrix-free construction of HSS representation\\
       using adaptive randomized sampling}

\author{Christopher Gorman\thanks{University of California, Santa Barbara,
    Santa Barbara, CA, USA. gorman@math.ucsb.edu} \and
  Gustavo Ch\'avez\thanks{Lawrence Berkeley National Laboratory,
    Berkeley, CA, USA. \{gichavez,pghysels,xsli\}@lbl.gov} \and
  Pieter Ghysels\footnotemark[2] \and
  Th{\'e}o Mary\thanks{University of Manchester, Manchester, UK.
    theo.mary@manchester.ac.uk} \and
  Fran\c{c}ois-Henry Rouet\thanks{Livermore Software Technology Corporation,
    Livermore, CA, USA. fhrouet@lstc.com} \and
  Xiaoye Sherry Li\footnotemark[2]}

\maketitle

\begin{abstract}
  We present new algorithms for randomized construction of
  hierarchically semi-separable matrices, addressing several practical
  issues. The HSS construction algorithms use a partially matrix-free,
  adaptive randomized projection scheme to determine the maximum
  off-diagonal block rank. We develop both relative and absolute stopping
  criteria to determine the minimum dimension of the random
  projection matrix that is sufficient for desired accuracy.
  Two strategies are discussed to adaptively
  enlarge the random sample matrix: repeated doubling of the number of
  random vectors, and iteratively incrementing the number of
  random vectors by a fixed number. The relative and absolute
  stopping criteria are based on probabilistic bounds for the
  Frobenius norm of the random projection of the Hankel blocks of the
  input matrix. We discuss parallel implementation and computation and
  communication cost of both variants. Parallel numerical results for
  a range of applications, including boundary element method matrices
  and quantum chemistry Toeplitz matrices, show the effectiveness,
  scalability and numerical robustness of the proposed algorithms.
\end{abstract}

\begin{keyword}
Hierarchically Semi-Separable,
randomized sampling,
randomized projection
\end{keyword}

%% file: tex/introduction.tex
%%%%%%%%%%%%%%%%%%%%%%%%%%%%%%%%%%%%%%%%%%%%%%%%%%%%%%%%%%%%%%%%%%%%%%%%
%%% Introduction
%%%%%%%%%%%%%%%%%%%%%%%%%%%%%%%%%%%%%%%%%%%%%%%%%%%%%%%%%%%%%%%%%%%%%%%%

\section{Introduction}

Randomization schemes have proven to be a powerful tool for computing a low-rank approximation of a dense matrix, or as we call it in this work, \textit{compressing} it. The main advantage of randomization is that these methods usually require fewer computations and communication than their traditional deterministic counterparts, resulting in large savings in terms of memory and floating point operations.

For classes of dense matrices that have off-diagonal blocks
that can be approximated as low-rank submatrices, randomized methods are particularly advantageous. 
These matrices are referred to as \emph{structured}, and there are many
types of matrix formats that can take advantage of this
structure; these include, to name a few, Hierarchically Semi-Separable (HSS) matrices~\cite{Chandrasekaran2005HSS}, $\mathcal{H}$ and $\mathcal{H}^{2}$ matrices~\cite{Hackbusch2000,Hackbusch2002}. This work focuses on HSS representations, and more specifically on efficient HSS compression. HSS compression is the central component of the HSS framework; once a matrix is compressed into its HSS form, one can take advantage of fast algorithms for multiplication, factorization, \ldots.

One way to speedup HSS compression involves using
randomization~\cite{martinsson2011fast,RandomReview2011}.
Randomization involves generating \emph{samples} of size
at least the maximum rank of the HSS representation.
Since the exact rank of low-rank matrices is usually not known in
practice, adaptive algorithms are needed in order to generate
sufficient, yet not too many, random samples, until the range is well approximated
and the matrix is compressed to a desired tolerance.  This ensures
robustness and high performance of the overall algorithm.

This paper builds on our previous work~\cite{FHRdistributedHSS},
which gives an explicit adaptive algorithm. One of the highlights of this work is the development of a new stopping criterion
that considers both relative and absolute error.
We demonstrate the effectiveness of this novel approach, and others, in a set of numerical experiments that showcase the scalability and robustness of the new algorithms on a variety of matrices from different applications.

The paper is organized as follows.
In Section~\ref{sec:HSS}, we discuss the HSS randomized construction algorithm, the two adaptive sampling strategies, and the new stopping criterion for terminating adaptation.
The parallel algorithms are presented and analyzed in
Section~\ref{sec:parallel}, followed by numerical experiments in
Section~\ref{sec:experiments}.
The probability theory proofs necessary for our new stopping
criterion are relegated to Appendix~\ref{sec:appendix}.

%% file: tex/hss.tex
%%%%%%%%%%%%%%%%%%%%%%%%%%%%%%%%%%%%%%%%%%%%%%%%%%%%%%%%%%%%%%%%%%%%%%%%
%%% HSS and Randomization

\section{Hierarchically Semi-Separable Matrices and Randomized Construction
\label{sec:HSS}}
In Section~\ref{ssec:HSSrepresentation}, we briefly describe the HSS
matrix format. Section~\ref{ssec:randomHSS} outlines the randomized
HSS construction scheme
from~\cite{martinsson2011fast}. Section~\ref{ssec:AdaptiveSampling}
introduces two new schemes that can be used to adaptively determine
the maximum HSS rank, making the randomized HSS construction more
robust. The final sections discuss the derivation of the stopping
criteria used in the adaptive schemes.

\subsection{HSS Representation\label{ssec:HSSrepresentation}}
This short HSS introduction mostly follows the notation
from~\cite{martinsson2011fast}. For a more detailed description of the
HSS matrix format, as well as for fast HSS algorithms, we refer the
reader to the standard HSS
references~\cite{Chandrasekaran2005HSS,Xia2010Superfast}.

The following notation is used: `$:$' is matlab-like notation for all
indices in the range, $^*$ denotes complex conjugation, $\#I_\tau$ is
the number of elements in index set
$I_{\tau} = \{ i_1, i_2, \cdots, i_n\}$, $R_\tau = R(I_\tau,:)$ is the
matrix consisting of only the rows $I_\tau$ of matrix $R$ and
$I_{\tau} \setminus I_{\nu}$ is the set of indices in $I_{\tau}$ minus
those in $I_{\nu}$.

Consider a square matrix $A \in \mathbb{C}^{N \times N}$ with an index
set $I_{A} = \{1,\dots,N \}$ associated with it. Let $\mathcal{T}$ be
a binary tree, ordered level by level, starting with zero at the root
node. Each node $\tau$ of the tree is associated with a contiguous
subset $I_\tau \subset \mathcal{I}$. For two siblings in the tree,
$\nu_1$ and $\nu_2$, children of $\tau$, it holds that
$I_{\nu_1} \cup I_{\nu_2} = I_\tau$ and
$I_{\nu_1} \cap I_{\nu_2} = \emptyset$. It follows that
$\cup_{\tau=\textrm{leaves}(\mathcal{T})} I_\tau =
I_{\textrm{root}(\mathcal{T})} = I_A$. The same tree $\mathcal{T}$ is
used for the rows and the columns of $A$ and only diagonal blocks are
partitioned. An example of the resulting matrix partitioning is given
in Figure~\ref{fig:HSS} and the corresponding tree is shown in
Figure~\ref{fig:HSStree}.

\begin{figure}
\centering
\begin{minipage}{.47\columnwidth}
  %\centering
  \begin{tikzpicture}[scale=4.0]
    \path[use as bounding box] (0,0) rectangle (1.2,1); % adjust to fit
    \node [label] at (0.25,0.25) {$A_{2,1}$};
    \node [label] at (0.75,0.75) {$A_{1,2}$};
    \draw (0,0) rectangle (1,1);
    \draw (0,0) rectangle (0.5,0.5);
    \draw (1,1) rectangle (0.5,0.5);
    \draw (0,1) rectangle (0.25,0.75);
    \draw [fill=gray] (0,1) rectangle (0.25/2,0.75+0.25/2);
    \draw [fill=gray] (0.25/2,0.75+0.25/2) rectangle (0.25,0.75);
    \draw (0.5,0.5) rectangle (0.75,0.25);
    \draw [fill=gray] (0.5,0.5) rectangle (0.75-0.25/2,0.25+0.25/2);
    \draw [fill=gray] (0.75-0.25/2,0.25+0.25/2) rectangle (0.75,0.25);
    \draw (0.25,0.75) rectangle (0.5,0.5);
    \draw [fill=gray] (0.25,0.75) rectangle (0.5-0.25/2,0.5+0.25/2);
    \draw [fill=gray] (0.5-0.25/2,0.5+0.25/2) rectangle (0.5,0.5);
    \draw (0.75,0.25) rectangle (1,0);
    \draw [fill=gray] (0.75,0.25) rectangle (0.75+0.25/2,0.25/2);
    \draw [fill=gray] (0.75+0.25/2,0.25/2) rectangle (1,0);
    \path [fill=gray] (0.02,0.35) rectangle (0.07,0.47);
    \path [fill=gray] (0.08,0.42) rectangle (0.13,0.47);
    \path [fill=gray] (0.14,0.42) rectangle (0.22,0.47);
    \path [fill=gray] (0.02+0.5,0.35+0.5) rectangle (0.07+0.5,0.47+0.5);
    \path [fill=gray] (0.08+0.5,0.42+0.5) rectangle (0.13+0.5,0.47+0.5);
    \path [fill=gray] (0.14+0.5,0.42+0.5) rectangle (0.22+0.5,0.47+0.5);
    \path [fill=gray] (0.25+0.01+0.5,0.01) rectangle (0.25+0.03+0.5,0.25/2-0.01);
    \path [fill=gray] (0.25+0.035+0.5,0.25/2-0.025-0.01) rectangle (0.25+0.055+0.5,0.25/2-0.01);
    \path [fill=gray] (0.25+0.06+0.5,0.25/2-0.025-0.01) rectangle (0.25+0.25/2-0.01+0.5,0.25/2-0.01);
    \path [fill=gray] (0.25/2+0.25+0.01+0.5,0.25/2+0.01) rectangle (0.25/2+0.25+0.03+0.5,0.25/2+0.25/2-0.01);
    \path [fill=gray] (0.25/2+0.25+0.035+0.5,0.25/2+0.25/2-0.025-0.01) rectangle (0.25/2+0.25+0.055+0.5,0.25/2+0.25/2-0.01);
    \path [fill=gray] (0.25/2+0.25+0.06+0.5,0.25/2+0.25/2-0.025-0.01) rectangle (0.25/2+0.25+0.25/2-0.01+0.5,0.25/2+0.25/2-0.01);
    \path [fill=gray] (0.01+0.5,0.01+0.25) rectangle (0.03+0.5,0.25/2-0.01+0.25);
    \path [fill=gray] (0.035+0.5,0.25/2-0.025-0.01+0.25) rectangle (0.055+0.5,0.25/2-0.01+0.25);
    \path [fill=gray] (0.06+0.5,0.25/2-0.025-0.01+0.25) rectangle (0.25/2-0.01+0.5,0.25/2-0.01+0.25);
    \path [fill=gray] (0.25/2+0.01+0.5,0.25/2+0.01+0.25) rectangle (0.25/2+0.03+0.5,0.25/2+0.25/2-0.01+0.25);
    \path [fill=gray] (0.25/2+0.035+0.5,0.25/2+0.25/2-0.025-0.01+0.25) rectangle (0.25/2+0.055+0.5,0.25/2+0.25/2-0.01+0.25);
    \path [fill=gray] (0.25/2+0.06+0.5,0.25/2+0.25/2-0.025-0.01+0.25) rectangle (0.25/2+0.25/2-0.01+0.5,0.25/2+0.25/2-0.01+0.25);
    \path [fill=gray] (-0.5+0.25+0.01+0.5,0.5+0.01) rectangle (-0.5+0.25+0.03+0.5,0.5+0.25/2-0.01);
    \path [fill=gray] (-0.5+0.25+0.035+0.5,0.5+0.25/2-0.025-0.01) rectangle (-0.5+0.25+0.055+0.5,0.5+0.25/2-0.01);
    \path [fill=gray] (-0.5+0.25+0.06+0.5,0.5+0.25/2-0.025-0.01) rectangle (-0.5+0.25+0.25/2-0.01+0.5,0.5+0.25/2-0.01);
    \path [fill=gray] (-0.5+0.25/2+0.25+0.01+0.5,0.5+0.25/2+0.01) rectangle (-0.5+0.25/2+0.25+0.03+0.5,0.5+0.25/2+0.25/2-0.01);
    \path [fill=gray] (-0.5+0.25/2+0.25+0.035+0.5,0.5+0.25/2+0.25/2-0.025-0.01) rectangle (-0.5+0.25/2+0.25+0.055+0.5,0.5+0.25/2+0.25/2-0.01);
    \path [fill=gray] (-0.5+0.25/2+0.25+0.06+0.5,0.5+0.25/2+0.25/2-0.025-0.01) rectangle (-0.5+0.25/2+0.25+0.25/2-0.01+0.5,0.5+0.25/2+0.25/2-0.01);
    \path [fill=gray] (-0.5+0.01+0.5,0.5+0.01+0.25) rectangle (-0.5+0.03+0.5,0.5+0.25/2-0.01+0.25);
    \path [fill=gray] (-0.5+0.035+0.5,0.5+0.25/2-0.025-0.01+0.25) rectangle (-0.5+0.055+0.5,0.5+0.25/2-0.01+0.25);
    \path [fill=gray] (-0.5+0.06+0.5,0.5+0.25/2-0.025-0.01+0.25) rectangle (-0.5+0.25/2-0.01+0.5,0.5+0.25/2-0.01+0.25);
    \path [fill=gray] (-0.5+0.25/2+0.01+0.5,0.5+0.25/2+0.01+0.25) rectangle (-0.5+0.25/2+0.03+0.5,0.5+0.25/2+0.25/2-0.01+0.25);
    \path [fill=gray] (-0.5+0.25/2+0.035+0.5,0.5+0.25/2+0.25/2-0.025-0.01+0.25) rectangle (-0.5+0.25/2+0.055+0.5,0.5+0.25/2+0.25/2-0.01+0.25);
    \path [fill=gray] (-0.5+0.25/2+0.06+0.5,0.5+0.25/2+0.25/2-0.025-0.01+0.25) rectangle (-0.5+0.25/2+0.25/2-0.01+0.5,0.5+0.25/2+0.25/2-0.01+0.25);
    \path [fill=gray] (0.015,0.15+0.5) rectangle (0.05,0.235+0.5);
    \path [fill=gray] (0.06,0.2+0.5) rectangle (0.09,0.235+0.5);
    \path [fill=gray] (0.1,0.2+0.5) rectangle (0.16,0.235+0.5);
    \path [fill=gray] (0.015+0.25,0.15+0.75) rectangle (0.05+0.25,0.235+0.5+0.25);
    \path [fill=gray] (0.06+0.25,0.2+0.75) rectangle (0.09+0.25,0.235+0.5+0.25);
    \path [fill=gray] (0.1+0.25,0.2+0.75) rectangle (0.16+0.25,0.235+0.5+0.25);
    \path [fill=gray] (0.5+0.015,0.15+0.5-0.5) rectangle (0.5+0.05,0.235+0.5-0.5);
    \path [fill=gray] (0.5+0.06,0.2+0.5-0.5) rectangle (0.5+0.09,0.235+0.5-0.5);
    \path [fill=gray] (0.5+0.1,0.2+0.5-0.5) rectangle (0.5+0.16,0.235+0.5-0.5);
    \path [fill=gray] (0.5+0.015+0.25,0.15+0.75-0.5) rectangle (0.5+0.05+0.25,0.235+0.5+0.25-0.5);
    \path [fill=gray] (0.5+0.06+0.25,0.2+0.75-0.5) rectangle (0.5+0.09+0.25,0.235+0.5+0.25-0.5);
    \path [fill=gray] (0.5+0.1+0.25,0.2+0.75-0.5) rectangle (0.5+0.16+0.25,0.235+0.5+0.25-0.5);
    \node [label] at (0.125/2*6,1-0.125/2*2) {\scriptsize $A_{3,4}$};
    \node [label] at (0.125/2*2,1-0.125/2*6) {\scriptsize $A_{4,3}$};
    \node [label] at (0.125/2*10,1-0.125/2*14) {\scriptsize $A_{6,5}$};
    \node [label] at (0.125/2*14,1-0.125/2*10) {\scriptsize $A_{5,6}$};
    \node [label] at (0.125/2,1-0.125/2) {\scriptsize $D_{7}$};
    \node [label] at (0.125/2*3,1-0.125/2*3) {\scriptsize $D_{8}$};
    \node [label] at (0.125/2*5,1-0.125/2*5) {\scriptsize $D_{9}$};
    \node [label] at (0.125/2*7,1-0.125/2*7) {\scriptsize $D_{10}$};
    \node [label] at (0.125/2*9,1-0.125/2*9) {\scriptsize $D_{11}$};
    \node [label] at (0.125/2*11,1-0.125/2*11) {\scriptsize $D_{12}$};
    \node [label] at (0.125/2*13,1-0.125/2*13) {\scriptsize $D_{13}$};
    \node [label] at (0.125/2*15,1-0.125/2*15) {\scriptsize $D_{14}$};
    \draw [decorate,decoration={brace,amplitude=2pt},yshift=0pt,xshift=-.1cm] (1.125,1-0.125*0) -- (1.125,1-0.125) node [black,midway,xshift=.3cm] {\footnotesize $I_7$};
    \draw [decorate,decoration={brace,amplitude=2pt},yshift=0pt,xshift=-.1cm] (1.125,1-0.125*1) -- (1.125,1-0.125*2) node [black,midway,xshift=.3cm] {\footnotesize $I_8$};
    \draw [decorate,decoration={brace,amplitude=2pt},yshift=0pt,xshift=-.1cm] (1.125,1-0.125*2) -- (1.125,1-0.125*3) node [black,midway,xshift=.3cm] {\footnotesize $I_9$};
    \draw [decorate,decoration={brace,amplitude=2pt},yshift=0pt,xshift=-.1cm] (1.125,1-0.125*3) -- (1.125,1-0.125*4) node [black,midway,xshift=.3cm] {\footnotesize $I_{10}$};
    \draw [decorate,decoration={brace,amplitude=2pt},yshift=0pt,xshift=-.1cm] (1.125,1-0.125*4) -- (1.125,1-0.125*5) node [black,midway,xshift=.3cm] {\footnotesize $I_{11}$};
    \draw [decorate,decoration={brace,amplitude=2pt},yshift=0pt,xshift=-.1cm] (1.125,1-0.125*5) -- (1.125,1-0.125*6) node [black,midway,xshift=.3cm] {\footnotesize $I_{12}$};
    \draw [decorate,decoration={brace,amplitude=2pt},yshift=0pt,xshift=-.1cm] (1.125,1-0.125*6) -- (1.125,1-0.125*7) node [black,midway,xshift=.3cm] {\footnotesize $I_{13}$};
    \draw [decorate,decoration={brace,amplitude=2pt},yshift=0pt,xshift=-.1cm] (1.125,1-0.125*7) -- (1.125,1-0.125*8) node [black,midway,xshift=.3cm] {\footnotesize $I_{14}$};

    \draw [decorate,decoration={brace,amplitude=3pt},yshift=0pt,xshift=-.1cm] (1.25,1-0.25*0) -- (1.25,1-0.25) node [black,midway,xshift=.3cm] {\footnotesize $I_3$};
    \draw [decorate,decoration={brace,amplitude=3pt},yshift=0pt,xshift=-.1cm] (1.25,1-0.25*1) -- (1.25,1-0.25*2) node [black,midway,xshift=.3cm] {\footnotesize $I_4$};
    \draw [decorate,decoration={brace,amplitude=3pt},yshift=0pt,xshift=-.1cm] (1.25,1-0.25*2) -- (1.25,1-0.25*3) node [black,midway,xshift=.3cm] {\footnotesize $I_5$};
    \draw [decorate,decoration={brace,amplitude=3pt},yshift=0pt,xshift=-.1cm] (1.25,1-0.25*3) -- (1.25,1-0.25*4) node [black,midway,xshift=.3cm] {\footnotesize $I_6$};

    \draw [decorate,decoration={brace,amplitude=4pt},yshift=0pt,xshift=-.1cm] (1.375,1-0.5*0) -- (1.375,1-0.5) node [black,midway,xshift=.3cm] {\footnotesize $I_1$};
    \draw [decorate,decoration={brace,amplitude=4pt},yshift=0pt,xshift=-.1cm] (1.375,1-0.5*1) -- (1.375,1-0.5*2) node [black,midway,xshift=.3cm] {\footnotesize $I_2$};

    \draw [decorate,decoration={brace,amplitude=4pt},yshift=0pt,xshift=-.1cm] (1.5,1) -- (1.5,0) node [black,midway,xshift=.3cm] {\footnotesize $I_0$};
  \end{tikzpicture}
  \caption{\footnotesize Illustration of an HSS matrix using $4$
    levels. Diagonal blocks are partitioned recursively. Gray blocks
    denote the basis matrices.}
  \label{fig:HSS}
\end{minipage} \hfill
\begin{minipage}{.47\columnwidth}
  \centering
  \begin{tikzpicture}[scale=.9]\scriptsize
    \node[cnode](0){0}
    child{node[cnode](1){1}
      child{node[cnode](3){3}   child{node[cnode](7){7}} child{node[cnode](8){8}} }
      child{node[cnode](4){4}   child{node[cnode](9){9}} child{node[cnode](10){10}} }
    }
    child{node[cnode](2){2}
      child{node[cnode](5){5}   child{node[cnode](11){11}}   child{node[cnode](12){12}} }
      child{node[cnode](6){6} child{node[cnode](13){13}} child{node[cnode](14){14}} }
    };
  \end{tikzpicture}
  \caption{\footnotesize Tree for Figure~\ref{fig:HSS}, using
    level-by-level top-down numbering. All nodes except the root store
    $U_\mu$ and $V_\mu$. Leaves store $D_\mu$, non-leaves
    $B_{\nu_1,\nu_2}$, $B_{\nu_2,\nu_1}$}
  \label{fig:HSStree}
\end{minipage}
\end{figure}

The diagonal blocks of $A$, denoted $D_\tau$, are stored as dense
matrices in the leaves $\tau$ of the tree $\mathcal{T}$
\begin{equation}
  D_{\tau} = A(I_{\tau},I_{\tau}) \, .
\end{equation}
The off-diagonal blocks $A_{\nu_1,\nu_2} = A(I_{\nu_1},I_{\nu_2})$,
where $\nu_1$ and $\nu_2$ denote two siblings in the tree, are
factored (approximately) as
\begin{equation}
  A_{\nu_1,\nu_2} \approx U^{\mathrm{big}}_{\nu_1} B_{\nu_1,\nu_2} \left( V_{\nu_2}^{\mathrm{big}}\right)^*  \, .
\end{equation}
The matrices $U^{\mathrm{big}}_{\nu_1}$ and
$V_{\nu_2}^{\mathrm{big}}$, which form bases for the column and row
spaces of $A_{\nu_1,\nu_2}$, are typically tall and skinny, with
$U^{\mathrm{big}}_{\nu_1}$ having $\#I_{\nu_1}$ rows and $r^r_{\nu_1}$
(column-rank) columns, $V_{\nu_2}^{\mathrm{big}}$ has $\#I_{\nu_2}$
rows and $r^c_{\nu_2}$ (row-rank) columns and hence $B_{\nu_1,\nu_2}$
is $r^r_{\nu_1} \times r^c_{\nu_2}$. The HSS-rank $r$ of matrix $A$ is
defined as the maximum of $r^r_\tau$ and $r^c_\tau$ over all
off-diagonal blocks, where typically $r \ll N$. The matrices
$B_{\nu_1,\nu_2}$ and $B_{\nu_2,\nu_1}$ are stored in the parent of
$\nu_1$ and $\nu_2$. For a non-leaf node $\tau$ with children $\nu_1$
and $\nu_2$, the basis matrices $U^{\mathrm{big}}_\tau$ and
$V^{\mathrm{big}}_\tau$ are not stored directly since they can be
represented hierarchically as
\begin{equation}\label{eq:nestedUV}
  U^{\mathrm{big}}_{\tau} = \begin{bmatrix} U^{\mathrm{big}}_{\nu_1} & 0 \\ 0 & U^{\mathrm{big}}_{\nu_2} \end{bmatrix} U_\tau
  \quad \textrm{and} \quad
  V^{\mathrm{big}}_{\tau} = \begin{bmatrix} V^{\mathrm{big}}_{\nu_1} & 0 \\ 0 & V^{\mathrm{big}}_{\nu_2} \end{bmatrix} V_\tau \, .
\end{equation}
Note that for a leaf node $U^{\mathrm{big}}_\tau = U_\tau$ and
$V^{\mathrm{big}}_\tau = V_\tau$. Hence, every node $\tau$ with
children $\nu_1$ and $\nu_2$, except for the root node, keeps matrices
$U_\tau$ and $V_{\tau}$. For example, the top two levels of the
example shown in Figure~\ref{fig:HSS} can be written out explicitly as
\begin{equation}
  A = \begin{bmatrix}
    D_3 & U_3 B_{3,4} V_4^* & \multicolumn{2}{c}{\multirow{2}{*}{$\begin{bmatrix} U_3 & 0 \\ 0 & U_4 \end{bmatrix} U_1 B_{1,2} V_2^* \begin{bmatrix} V_5^* & 0 \\ 0 & V_6^* \end{bmatrix}$}} \\
    U_4 B_{4,3} V_3^* & D_4 &  \\
    \multicolumn{2}{c}{\multirow{2}{*}{$\begin{bmatrix} U_5 & 0 \\ 0 & U_6 \end{bmatrix} U_2 B_{2,1} V_1^* \begin{bmatrix} V_3^* & 0 \\ 0 & V_4^* \end{bmatrix}$}}
    & D_5 & U_5 B_{5,6} V_6^* \\
      & & U_6 B_{6,5} V_5^* & D_6
    \end{bmatrix} \, .
    \label{eq:HSSexamplematrix}
\end{equation}
The storage requirement for an HSS matrix is
$\mathcal{O}(rN)$. Construction of the HSS generators will be
discussed in the next section. Once an HSS representation of a matrix
is available, it can be used to perform matrix-vector multiplication
in $\mathcal{O}(rN)$ operations compared to $\mathcal{O}(N^2)$ for
classical dense matrix-vector multiplication,
see~\cite{martinsson2011fast,FHRdistributedHSS}.

\subsection{Randomized HSS Construction\label{ssec:randomHSS}}
Here we review the randomized HSS construction algorithm as presented
by Martinsson~\cite{martinsson2011fast} as well as the adaptive
variant as described
in~\cite{FHRdistributedHSS,ghysels2014multifrontal,ghysels2017robust}. This
algorithm requires matrix vector products, for the random sampling, as
well as access to certain entries of the input matrix, for
construction of the diagonal blocks $D_i$ and the $B_{ij}$
generators. We therefore refer to it as ``partially matrix-free''. A
fully matrix-free randomized HSS compression algorithm was presented
in~\cite{Liu2016PRMatfree}, which does not require access to
individual matrix elements, but requires $\mathcal{O}(\log N)$ more
random projection vectors instead of requiring access to
$\mathcal{O}(rN)$ matrix elements.
% \todo{quantify how many more?, also not adaptive!}  \chgcomment{See
% Liu 2016, Section 5.1 for computational cost}

Let us assume for now that the maximum rank encountered in any of the
HSS off-diagonal blocks is known to be $r$. Let $R$ be an $N \times d$
random matrix, where $d = r + p$ with $p$ a small oversampling
parameter. The matrix $U_\tau$ in the HSS representation is a basis
for the column space of the off-diagonal row block
$A(I_\tau, I_A \setminus I_\tau)$. Likewise, $V_\tau$ is a basis for
the row space of $A(I_A \setminus I_\tau, I_\tau)^*$. These row and
column off-diagonal blocks are called Hankel blocks. In order to
compute the random projection of these Hankel blocks, we first compute
both $S^r = A R$ and $S^c = A^* R$, such that $S^r(I_\tau, :)$ is the
random projection of the entire row block $A(I_\tau, :)$. Let $D_\tau$
for non-leaf nodes be defined (recursively) as
$D_\tau = \begin{bmatrix} D_{\nu_1} & A_{\nu_1,\nu_2} \\
  A_{\nu_2,\nu_1} & D_{\nu_2} \end{bmatrix}$, and with nodes
$\tau_i, \ldots, \tau_j$ on level $\ell$ of the HSS tree, we
define the block diagonal matrix
$D^{(\ell)} = \text{diag}\left( D_{\tau_i}, \dots, D_{\tau_j}
\right)$.  At each level of the HSS tree we can compute the samples of
the Hankel blocks as
\begin{eqnarray}
  S^{r,(\ell)} = \left( A - D^{(\ell)} \right) R = S^r - D^{(\ell)} R \quad \text{and} \quad
  S^{c,(\ell)} = S^c - D^{(\ell)*} R \, . \label{eq:sample_off_diag}
\end{eqnarray}
By starting from the leaf level of the HSS tree and working up towards
the root, this calculation can be performed efficiently, since the
off-diagonal blocks of $D_{\tau}$ are already compressed.  Consider a
leaf node $\tau$, and let
$S^r_\tau = S^r(I_\tau, :) - D_\tau R(I_\tau, :)$.  To compute the
$U_\tau$ basis from this random projection matrix, an interpolative
decomposition (ID) is used, which expresses the matrix $S^r_\tau$ as a
linear combination of a selected set $J_{\tau}^r$ of its rows
$S^r_\tau = U_{\tau} S^r_\tau(J^r_{\tau}, :) +
\mathcal{O(\varepsilon)}$, where $\#J^r_{\tau}$ is the
$\varepsilon$-rank of $S^r_{\tau}$. This decomposition can be computed
using a rank-revealing QR factorization, or QR with column pivoting
(QRCP), applied to $(S^r_\tau)^*$ as follows: % \todo{do this
% on the transpose, then define $U$ from this}
\begin{align} \label{eq:ID}
  (S^r_\tau)^* & \approx Q \begin{bmatrix} R_1 & R_2 \end{bmatrix} \Pi^T
                                                 = \left( Q R_1 \right) \begin{bmatrix} I & R_1^{-1} R_2 \end{bmatrix} \Pi^T \\
               &\approx (S^r_\tau(J^r_{\tau}, :))^* \begin{bmatrix} I & R_1^{-1} R_2 \end{bmatrix} \Pi^T
                                                                  \, , \label{eq:QRCPSr}
\end{align}
with $\Pi$ a permutation matrix moving columns $J^r_{\tau}$ to the
front, $Q$ orthogonal and $R_1$ upper triangular. Note that QRCP
selects columns $J^r_{\tau}$ of $(S^r_\tau)^*$, which corresponds to
the transposed rows $J^r_{\tau}$ of
$S^r_\tau$. From~(\ref{eq:QRCPSr}), $U_\tau$ can be defined as
\begin{equation}
  U_\tau = \Pi \begin{bmatrix} I \\ \left( R_1^{-1} R_2 \right)^* \end{bmatrix}
  = \Pi^r_\tau \begin{bmatrix} I \\ E^r_\tau \end{bmatrix} \, ,
\end{equation}
and, likewise
$V_\tau = \Pi^c_{\tau} \begin{bmatrix} I &
  (E^c_\tau)^* \end{bmatrix}^*$. From these definitions of $U_\tau$
and $V_\tau$, we see that if $\nu_1$ and $\nu_2$ are leaf nodes
\begin{equation}
  A_{\nu_1,\nu_2} \approx U_{\nu_1} B_{\nu_1,\nu_2} \left( V_{\nu_2} \right)^* =
  \Pi^r_\tau \begin{bmatrix} I \\ E^r_\tau \end{bmatrix}
  B_{\nu_1,\nu_2} \begin{bmatrix} I & (E^c_\tau)^* \end{bmatrix} (\Pi^c_{\tau})^T
\end{equation}
and hence $B_{\nu_1,\nu_2} = A_{\nu_1,\nu_2}(J^r_{\nu_1},J^c_{\nu_2})$
is a sub-block of $A_{\nu_1,\nu_2}$.

Note that Equation~(\ref{eq:sample_off_diag}) can be used for the leaf
levels. On the higher levels however, we need to guarantee the nested
basis property, see Equation~(\ref{eq:nestedUV}). For a non-leaf node
$\tau$ with (leaf) children $\nu_1$ and $\nu_2$ we have
\begin{align}
  \label{eq:Sr_start}
  S^{r,(\ell)}_{\tau} &= \left( A(I_{\tau}, :) - D_{\tau} \right) R
  % = A(I_{\tau}, :) R - D_{\tau} R(I_{\tau},:) \\
  = A(I_{\tau}, :) R
  - \begin{bmatrix} D_{\nu_1} & A_{\nu_1,\nu_2}\\ A_{\nu_2,\nu_1} & D_{\nu_2}\end{bmatrix}
                                                                    \begin{bmatrix} R(I_{\nu_1},:) \\ R(I_{\nu_2},:) \end{bmatrix} \\
                      &= \begin{bmatrix} S^r_{\nu_1} \\ S^{r}_{\nu_2} \end{bmatrix}
  - \begin{bmatrix}  & A_{\nu_1,\nu_2}\\ A_{\nu_2,\nu_1} & \end{bmatrix}
                                                           \begin{bmatrix} R(I_{\nu_1},:) \\ R(I_{\nu_2},:) \end{bmatrix} \\
                                                           % = \begin{bmatrix} S^r_{\nu_1} - A_{\nu_1,\nu_2} R(I_{\nu_2},:) \\
                                                           %   S^r_{\nu_2} - A_{\nu_2,\nu_1} R(I_{\nu_1},:) \end{bmatrix} \\
  \label{eq:Sr_approx} &\approx \begin{bmatrix} U_{\nu_1} S^r_{\nu_1}(J^r_{\nu_1},:) - U_{\nu_1} B_{\nu_1,\nu_2} V_{\nu_1}^* R(I_{\nu_2},:) \\
                        U_{\nu_2} S^r_{\nu_2}(J^r_{\nu_1},:) - U_{\nu_2} B_{\nu_2,\nu_1} V_{\nu_1}^* R(I_{\nu_1},:) \end{bmatrix} \\
                      &\approx \begin{bmatrix} U_{\nu_1} \\ & U_{\nu_2} \end{bmatrix}
                                                        \begin{bmatrix} S^r_{\nu_1}(J^r_{\nu_1},:) - B_{\nu_1,\nu_2} V_{\nu_2}^* R(I_{\nu_2},:) \\
                                                          S^r_{\nu_2}(J^r_{\nu_1},:) - B_{\nu_2,\nu_1} V_{\nu_1}^* R(I_{\nu_1},:) \end{bmatrix} \, .
  \label{eq:U1U2Sr}
\end{align}
We let $\tau.R^r \gets V_\tau^* R(I_{\tau},:)$ and
\begin{equation} \label{eq:definition_of_tauSr}
  \tau.S^r \gets \begin{bmatrix} S^r_{\nu_1}(J^r_{\nu_1},:) - B_{\nu_1,\nu_2} \nu_2.R \\
    S^r_{\nu_2}(J^r_{\nu_1},:) - B_{\nu_2,\nu_1} \nu_1.R \end{bmatrix} \, ,
\end{equation}
so we can apply $\text{ID}(\tau.S^r)$ in order to compute $U_\tau$,
since the $U_{\nu_1}$ and $U_{\nu_2}$ bases have been factored out in
Eq.~(\ref{eq:U1U2Sr}). This can be applied recursively to nodes with
children that are non-leaf nodes. A similar reasoning also applies for
$\text{ID}(\tau.S^c)$.

\input{tex/tables/defs.tex}

\input{tex/algo/hss_comp.tex}

\input{tex/algo/comp_local_reduce_2.tex}

These steps, Eqs.~(\ref{eq:sample_off_diag})-(\ref{eq:U1U2Sr}), are
implemented in Algorithms~\ref{algo:HSScompress}
and~\ref{algo::computelocalsamples}, with some definitions in
Table~\ref{tab::addfunctions}. We use notation like $\tau.J^r$ to
denote that the temporary variable $J^r_{\tau}$ is stored at node
$\tau$. Algorithm~\ref{algo:HSScompress} implements the recursive
bottom-up HSS tree traversal, calling the helper functions
\texttt{ComputeLocalSamples} (computation of $\tau.S^r$ and
$\tau.S^c$, Eqs.~(\ref{eq:sample_off_diag}) and~(\ref{eq:U1U2Sr})) and
\texttt{ReduceLocalSamples} (compution of $\tau.R^r$ and $\tau.R^c$)
which are defined in Algorithm~\ref{algo::computelocalsamples}.  After
successful HSS compression, each non-leaf HSS node $\tau$ stores
$\tau.B_{\nu_1,\nu_2}$, and each non-leaf node stores $\tau.U$ and
$\tau.V$, although in practice, $\tau.U$ is stored using $\tau.\Pi^r$
and $\tau.E^r$.

\input{tex/algo/adaptive_comp_node.tex}

\begin{figure}[H]
\centering
\includegraphics[width=0.7\linewidth]{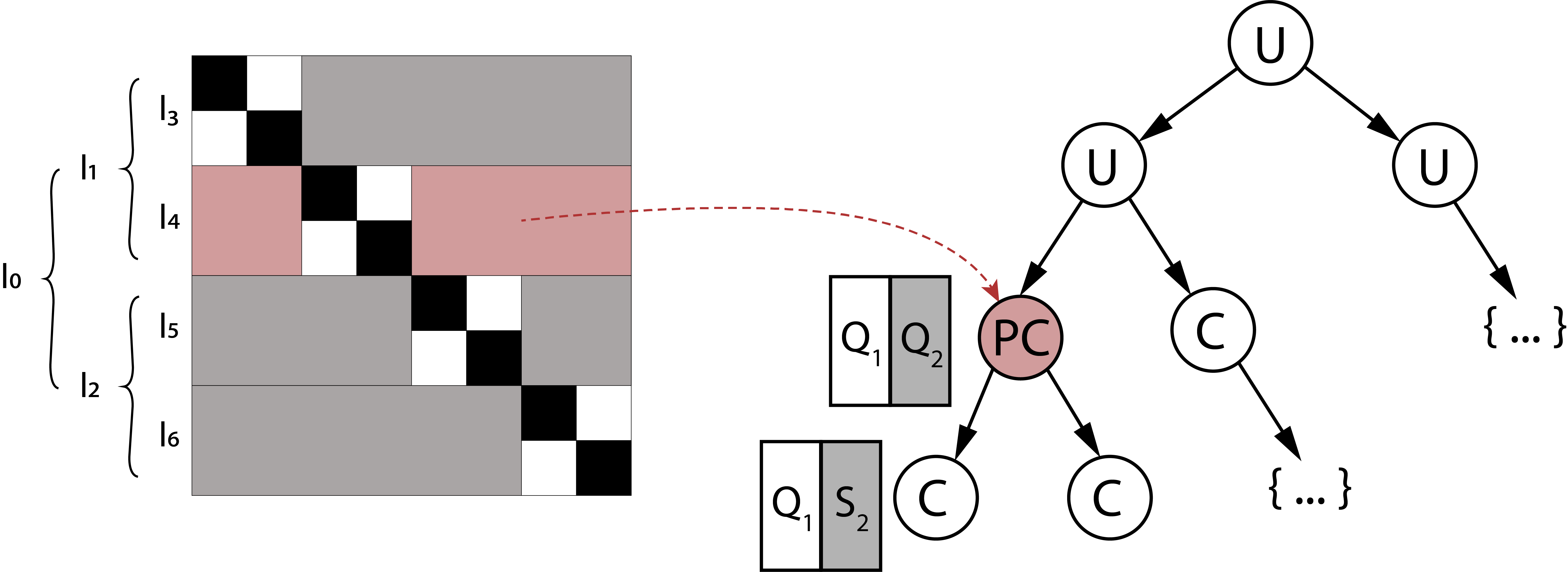}
\caption{Illustration of the tree traversal during adaptive HSS
  compression. Each node in the tree can be in either
  \texttt{UNTOUCHED} (U), {\PC} (PC), or \texttt{COMPRESSED} (C)
  state. For the meaning of the $Q$ and $S$ matrices, see
  Section~\ref{ssec::incremental_adaptive}.}
\label{fig:compression stages}
\end{figure}

In~\cite{FHRdistributedHSS}, we have extended the randomized HSS
construction algorithm to make it adaptive, relaxing the condition
that the maximum HSS is known {\it a priori}. This adaptive scheme is
illustrated in Algorithm~\ref{algo::HSScompressNodeAdaptive} and it
works as follows. Each HSS node has a state field which can be
\texttt{UNTOUCHED}, \texttt{COMPRESSED} or {\PC}. Each node starts in
the \texttt{UNTOUCHED} state. The compression proceeds bottom-up from
the leaf nodes, as in the non-adaptive
Algorithm~\ref{algo:HSScompress}, with an initial number of random
vectors $d_0$. At each (non-root) node $\tau$, bases $U_\tau$ and
$V_\tau$ are computed using the ID, which, however, might fail if it
is detected that random projection with $d_0$ random vectors is not
sufficient to accurately capture the range of the corresponding Hankel
block with a prescribed relative or absolute tolerance
$\varepsilon_{\text{rel}}$ or $\varepsilon_{\text{abs}}$
respectively. How this can be detected will be discussed in detail in
the following sub-sections. If at a node in the HSS tree, it is
decided that $d_0$ is not sufficient, that node is marked as {\PC} and
the compression algorithm returns to the outer loop. Note that in a
parallel run, processing of other HSS nodes in independent sub-trees
can still continue. The number of random projection vectors is
increased by $\Delta d$ and the recursive HSS compression is called
again. At this point, there will be at least one node in the HSS tree
which is in the \texttt{PARTIALLY\_COMPRESSED} state. All of the
descendants of this node are in the \texttt{COMPRESSED} state (or
compression of that {\PC} could not have been started), and all of the
ancestors are in the \texttt{UNTOUCHED} state. The compression
algorithm will again start at the leaf nodes. For the already
\texttt{COMPRESSED} nodes, the $U$ and $V$ bases do not need to be
recomputed, but these nodes still need to be visited in order to add
$\Delta d$ extra columns to $\tau.S^r$/$\tau.S^c$
(Eqs.~(\ref{eq:sample_off_diag}) and~(\ref{eq:U1U2Sr})) and
$\tau.R^r$/$\tau.R^c$. For a graphical representation of the adaptive
compression procedure see Figure~\ref{fig:compression stages}.

Theorem 4.2 in~\cite{vogel2016superfast} shows that for a matrix $A$
and its HSS approximation $H$, with a total of $L$ levels, and with
each off-diagonal block compressed with tolerance $\varepsilon$, it
holds that $\| A - H \| \le L \varepsilon$.  Inspired by this result,
in the HSS compression, Algorithms~\ref{algo:HSScompress}
and~\ref{algo::HSScompressNodeAdaptive}, the user requested absolute
and relative tolerances $\varepsilon_a$ and $\varepsilon_r$ are scaled
with the current HSS level for the rank-revealing QR
factorization. Hence, it is expected that the final HSS compression
satisfies $\| A - H \| \le \varepsilon_a$ or
$\| A - H \|/\|A\| \le \varepsilon_r$.

\subsection{Difficulties with Adaptive Sampling}
In this section we discuss some difficulties associated with the
adaptive randomized HSS construction algorithm as described in
Section~\ref{ssec:randomHSS},
Algorithm~\ref{algo::HSScompressNodeAdaptive}. The main issue is to
decide whether the sample matrices $S^r$ and $S^c$ capture the column
space of the corresponding Hankel blocks well enough. As we shall see
later, this decision relies on a good estimate of the matrix norm and
hence
% \todo{norm of the Hankel block?} of the sample matrix,
the error estimation of the residual matrix (see
Eq.~\eqref{eq:stopping_criteria}).
% \todo{just refer to eq later on?}).
A good error estimation is critical in devising the stopping criteria
of adaptation.  The goal is to ensure sufficient samples are used to
guarantee the desired accuracy is, while not to perform too much
oversampling, as this degrades performance. We first review a
well-known result from the literature, and explain how it falls short
when used as a stopping criterion. Then in
Section~\ref{ssec:sampling_doubling} we propose our new approach.

From~\cite{RandomReview2011}, we have the following lemma:
\begin{lemma}[Lemma 4.1 in~\cite{RandomReview2011}]
\label{lem:RandomLem}
Let $B$ be a real $m\times n$ matrix. Fix a positive integer $p$ and a
real number $\alpha>1$. Draw an independent family
$\braces{\omega^{i}:i = 1,\cdots,p}$ of standard Gaussian
vectors. Then
\begin{equation}\label{eqn:HMT_err_bound}
  \norm{B}_{2} \le \alpha \sqrt{\frac{2}{\pi}} \max_{i=1,\ldots,p}
  \norm{B\omega^{i}}_{2}
\end{equation}
except with probability $\alpha^{-p}$.
\end{lemma}
We will refer to the parameter $p$ as the oversampling parameter. The
value of $\alpha$ can be used as a trade-off, making the bound less
strict or less probable. Consider for instance $\alpha = 10$ and
$p = 10$. Then the bound in~(\ref{eqn:HMT_err_bound}) holds with very
high probability $10^{-10}$. In~\cite{Liu2016PRMatfree},
relation~(\ref{eqn:HMT_err_bound}) is used as a stopping criterion in
an adaptive HSS construction algorithm, with $B = (I - QQ^{*})F$, as
\begin{equation} \label{eqn::HMT_IQQH}
  \norm{\parens{I - QQ^{*}}F}_{2} \le \alpha\sqrt{\frac{2}{\pi}}
  \max_{i=1,\ldots,p} \norm{\parens{I - QQ^{*}}F\omega^{i}}_{2}
  \le \eps
\end{equation}
where $F$ is a Hankel block of the input matrix $A$, $Q$ is a matrix
with orthonormal columns which approximates the range of $F$, and
$\eps$ is tolerance. The column dimension of $Q$ is increased until
the second inequality in~(\ref{eqn::HMT_IQQH}) is satisfied, which
guarantees that $\|(I-QQ^*)F\|_2 \le \eps$. There are two drawbacks
with this criterion: 1) Since~(\ref{eqn:HMT_err_bound}) gives an upper
bound, the error might be drastically overestimated; 2) if different
$F$ blocks vary greatly in size, then different blocks will be
compressed to different \emph{relative} tolerances, even though they
are compressed to the same \emph{absolute} tolerance.
%\todo{\chgcomment{If possible, it would be good to have more references
%    where this stopping criterion is used. It probably is used a lot
%    so that our critique is noted.}}

The potential overestimation of the error is not the main
difficulty. The error estimate in~\eqref{eqn::HMT_IQQH} assumes that
we have access to the entire matrix Hankel block $F$ and the
associated random samples. However, in the randomized matrix-free HSS
construction as described in Section~\ref{ssec:randomHSS}, only the
local random samples $\tau.S^r$ and $\tau.S^c$ are
available. In~\eqref{eq:Sr_start}-\eqref{eq:U1U2Sr}, the random sample
matrix $S_{\tau}^r$ of the Hankel block corresponding to HSS node
$\tau$, i.e., $A(I_{\tau}, I_A \setminus I_{\tau})$, is expressed in
terms of the basis generators $U_{\nu_1}$ and $U_{\nu_2}$ of its
children, see Eq.~\eqref{eq:U1U2Sr}. The random sample $S_{\tau}^r$ of
the Hankel block is constructed by subtracting the random sample
$D_{\tau} R(I_{\tau},:)$ of the diagonal block,
see~\eqref{eq:Sr_start}, using the already compressed
$D_{\tau}$ ($D_{\tau}$ is exact only at the leaves). Hence, this introduces an approximation error,
see~\eqref{eq:Sr_approx}. Another issue arises due to the HSS nested
basis property, see Eq.~\eqref{eq:nestedUV}. In order to maintain the
HSS nested basis property, the Hankel block is not compressed
directly, but only its coefficient in the
$[ U_{\nu_1} \,\, 0 \,; \,\, 0 \,\, U_{\nu_2}]$ basis is
compressed. This coefficient is defined
in~\eqref{eq:definition_of_tauSr} as $\tau.S^r$. The difficulty arises
from the use of the interpolative decomposition to compress the
off-diagonal blocks, see~\eqref{eq:ID}, since it leads to
non-orthonormal bases $U_{\tau}$. Hence, the bound from
Lemma~\ref{lem:RandomLem} should really be applied to
\begin{equation}
  \norm{
    \begin{bmatrix} U_{\nu_{1}} & 0 \\ 0 & U_{\nu_{2}} \end{bmatrix}
    \parens{I - Q_{\tau}Q_{\tau}^{*}}S_{\tau}}_{2} \, ,
  \label{eq:exact_off_diag}
\end{equation}
to derive a correct stopping criterion for adaptive HSS compression.
As long as the $U$ bases are orthonormal, there is no problem, since
$\norm{\cdot}_{2}$ is unitarily invariant.
% However, when we use ID as in Eq.~\eqref{eq:ID}, the $U$ bases are
% \emph{not} orthonormal.  \todo{We might be able to argue that U is
% ``almost'' orthonormal?  From the properties of the ID, namely that
% the U block is bounded, at least when using a proper (strong) rank
% revealing QR?}  \todo{\chgcomment{Yes, when using a (strong) rank
% revealing QR, the $U$ matrices will be well-conditioned and
% essentially orrthonormal, so it can probably be ignored. The
% absolute stopping criterion could be scaled in a way related to the
% dimension.  A relative stopping criterion may be better in practice
% because the non-orthonormality in the numerator and denominator
% ``cancel''.}}  This has caused theory and practice to deviate except
% at the lowest level of the HSS tree.
In practice, a (strong) rank-revealing QR
factorization~\cite{gu1996efficient} will cause the $U$ matrices to be
well-conditioned, and its elements to be
bounded~\cite{RandomReview2011}, so the absolute tolerance should not
be affected much.  It seems plausible for the non-orthonormal factors
to essentially cancel out when using the relative stopping criterion,
making it a more reliable estimate. Due to the hierarchical nature of
the $U$ and $V$ matrices, it would be possible to compute the matrix
products in Eq.~\eqref{eq:exact_off_diag} if this would be desired.

The error bound in Lemma~\ref{lem:RandomLem} is not conducive to
relative error bounds because it only bounds the error, not the matrix
norm. It is not possible to use this bound in an attempt to compute a
relative stopping criterion. In
Section~\ref{ssec::incremental_adaptive}, we propose a new stopping
criterion for adaptive HSS compression based on both absolute and
relative bounds, inspired by Lemma~\ref{lem:RandomLem}, and the
derivation of the accurate relative bound is given in
Section~\ref{ssec:math_theory}.

% \chgcomment{We probably do not need all of the above calculations,
% but I wanted to include what I thought would be the maximum needed,
% as it easier to delete some of the calculations.}

\subsection{Adaptive Sampling Strategies\label{ssec:AdaptiveSampling}}
While traversing the HSS tree during compression, a node may fail to
compress to the desired accuracy due to a lack of random samples; this
node, which corresponds to one off-diagonal Hankel block of the matrix
$A$, is then labeled as {\PC}. We consider two strategies to add more
random vectors: {\dbl}, and {\inc}.  There are advantages and
disadvantages to both approaches.
Algorithms~\ref{algo::adaptive_rs_dbl} ({\dbl})
and~\ref{algo::adaptive_rs_inc} ({\inc}) present the two strategies
for the computation of an approximate basis for the range of a matrix
$A$. Algorithms~\ref{algo::adaptive_rs_dbl}
and~\ref{algo::adaptive_rs_inc} correspond to the successive
(rank-revealing) QR factorizations at the {\PC} node, and is part of
the bigger Algorithm~\ref{algo::HSScompressNodeAdaptive}. The
successive random sampling steps are actually performed in
Algorithms~\ref{algo::HSScompressNodeAdaptive}, but are included here
for convenience. In this section,
Algorithms~\ref{algo::adaptive_rs_dbl} and~\ref{algo::adaptive_rs_inc}
are presented as only returning the orthonormal basis of the range of
$A$, but when used in the overall HSS compression,
Algorithm~\ref{algo::HSScompressNodeAdaptive}, also the upper
triangular factor from the final rank-revealing QR factorization is
required.

% \todo{better explain how these fit into the larger algo, use sizes of
%   the Hankel block, not of the reduces size local samples!!}

% The flop counts do not include the sampling flops.

% \xslcomment{It's worth to present them separately to explain the
% stopping criteria and the difference between the two strategies.  }

%% \xslcomment{The {\dbl} strategy described in the flop count
%% paragraph is different from my understanding. We should discuss.}

% \xslcomment{All the algorithm pseudocodes use {\inc} strategy.}
% \xslcomment{Shall we put flop count in this subsection?}

\subsubsection{{\dbl} Strategy with Oversampling\label{ssec:sampling_doubling}}

\input{tex/algo/adaptive-rs-dbl}

Algorithm~\ref{algo::adaptive_rs_dbl} computes an (orthogonal)
approximate basis for the range of an $m \times n$ matrix $A$, with an
oversampling parameter $p$, up to a given relative or absolute
tolerance $\varepsilon_r$ or $\varepsilon_a$. Initially, a random
projection of the input matrix $A$ with a tall and skinny random
matrix $R$, with $d_0 + p$ columns, is computed as $S = A R$.  Let $Q$
be an orthonormal basis for $S$. From~\cite[Theorem
1.1]{RandomReview2011}, we have
\begin{equation}\label{eq:expIQQtA}
  \mathbb{E}\|(I - QQ^*)A \| \le
  \left[1 + 4 \frac{\sqrt{d_0 + p}}{p-1} \cdot \sqrt{\min\{m, n\}} \right] \sigma_{d_0+1} \, .
\end{equation}
where $\sigma_{d_0+1}$ is the $(d_0+1)$-th singular value of $A$. The
factor in brackets in Eq.~\eqref{eq:expIQQtA}, the deviation from the
optimal error $\sigma_{d_0+1}$, decreases with increasing oversampling
$p$. In order to guarantee a relative or absolute error bound on
$\|A-QQ^*A\|$, we apply a modified RRQR factorization (a column
pivoted QR) to $S$, which includes this factor. The modified RRQR
(\texttt{RRQR\_HMT} in Algorithm~\ref{algo::adaptive_rs_dbl}) return a
rank $r=k$ as soon as
\begin{equation}\label{eq:RRQRHMTstop}
  R_{k,k} \le \frac{\varepsilon_a}
  {1 + 4 \frac{\sqrt{d_0 + p}}{d_0 + p - k -1} \cdot \sqrt{\min\{m, n\}} }
  \quad \text{or} \quad
  \frac{R_{k,k}}{R_{1,1}} \le \frac{\varepsilon_r}
  {1 + 4 \frac{\sqrt{d_0 + p}}{d_0 + p - k -1} \cdot \sqrt{\min\{m, n\}} } \, .
\end{equation}
This modification is to take into account that at step $k$ of the
\texttt{RRQR\_HMT} algorithm, the amount of oversampling is
$d_0+p - k - 1$. If the rank $r$ returned by the modified RRQR is
$r < d_0$, i.e., at least $p$ oversampling vectors are used, then the
orthonormal basis $Q$ returned by RRQR is accepted. However, if RRQR
does not achieve the required tolerance, or if it achieves the
required tolerance but $d_0 \leq r \leq d_0 + p$, then the resulting
$Q$ basis is rejected. The number of random vectors (excluding the $p$
oversampling columns) is doubled and the random projection with these
new vectors is added to $S$. A rank-revealing QR factorization is
applied to the entire new random projection matrix $S$. Since in this
scenario RRQR is always redone from scratch, the number of adaptation
steps should be minimized in order to minimize the work in
RRQR. Doubling the number of random vectors in each step ensures one
only needs $\mathcal{O}(\log(r))$ adaptation steps\footnote{Logarithm is in base
  $2$ throughout the paper.}, but in the
worst-case scenario, the amount of oversampling is $r / 2 - 1 + p$.

This is the approach used
in~\cite{FHRdistributedHSS,ghysels2014multifrontal,ghysels2017robust}
for adaptive randomized HSS construction, although in these
references, the scaling factor from Eq.~\eqref{eq:expIQQtA} was not
included. Including this scaling factor slightly increases the rank,
but gives more accurate compression. When used in the larger HSS
compression algorithm, the matrix sizes $m$ and $n$ from
Eq.~\eqref{eq:RRQRHMTstop} are the sizes of the original Hankel block.

\subsubsection{{\inc} Strategy based on Lemma~\ref{lem:RandomLem}}
\label{ssec::incremental_adaptive}

\input{tex/algo/adaptive-rs-inc}

Algorithm~\ref{algo::adaptive_rs_inc} presents a new approach to
adaptive rank determination. The adaptive scheme used in
Algorithm~\ref{algo::adaptive_rs_dbl} performs a rank-revealing
factorization in each adaptation step with a random projection matrix
that is double the size of that in the previous step. The
rank-revealing factorization has to be performed from scratch in each
step, which leads to additional computational cost and
communication. In particular, the communication involved in column
pivoting in the repeated RRQR factorizations can be a serious
bottleneck in a distributed memory code.
% \todo{in section 5, discuss the performance of RRQR, PDGEQPF,
% refer to alternatives}
% \todo[inline]{Describe this algorithm in detail; how to decide when
% to accept the compression! Block version of the algorithm from
% HMT. Why incrementing makes sense here, instead of doubling?}
In contrast, in Algorithm~\ref{algo::adaptive_rs_inc}, the
rank-revealing factorization is only performed when the number of
sample vectors is guaranteed to be sufficient. This decision criterion
is based on relation~\eqref{eqn::HMT_IQQH} and uses an adaptive
blocked implementation.

%\todo[inline]{Move this to the introduction?}
In~\cite{RandomReview2011}, Algorithm~4.2 (``Adaptive Randomized Range
Finder'') is presented for the adaptive computation of an orthonormal
basis $Q$ for the range of a matrix $A$, up to an absolute tolerance
$\varepsilon_a$. This algorithm computes an approximate orthonormal
basis for the range of $A$ one vector at a time and then determines
how well this basis approximates the range. Adding one vector at a
time to the basis amounts to performing multiple matrix-vector
products (BLAS-2), which are memory-bound and less efficient on modern
hardware. In~\cite{Liu2016PRMatfree}, Algorithm~2 (``Parallel adaptive
randomized orthogonalization'') uses an approach similar to
Algorithm~4.2 in~\cite{RandomReview2011}, but implements a blocked
version, relying on BLAS-3 operations which can achieve much higher
performance. This also relies on an absolute tolerance only.
Algorithm~\ref{algo::adaptive_rs_inc} implements a blocked version,
but it uses both an absolute and a relative stopping criterion. The
stopping criterion is based on a stochastic error bound and block
Gram-Schmidt orthogonalization~\cite{stewart2008,bjorck1994GS}.
% Our stopping criterion is described in detail
% in~\ref{algo::adaptive_rs_inc}. We highlight the main features here.

Let $Q$ be an orthonormal approximate basis for the range of an
$m \times n$ matrix $A$. Given a random Gaussian matrix
$R \in \mathbb{R}^{n\times d}$, compute first the random samples
$S = AR$ followed by the projection
$\widehat{S} = \parens{I - QQ^{*}}S$. This last operation is a step of
block Gram-Schmidt. In practice, to ensure orthogonality, we apply
this step twice~\cite{stewart2008}. Thus, $\widehat{S}$ contains
information about the range of $A$ that is not included in $Q$.  A
small $||\widehat{S}||_{F}$ means that either $Q$ was already a good
basis for the range of $A$, or $S$ was in the range of $Q$, which is
unlikely. If $||\widehat{S}||_{F}$ is not small enough, $\widehat{S}$
is used to expand $Q$. Since $\widehat{S}$ is already orthogonal to
$Q$, $\widehat{S}$ only needs to be orthogonalized using a QR
factorization, $[\overline{Q}, R] = \text{QR}(\widehat{S})$, and can
then be added to $Q \gets [Q \,\, \overline{Q}]$.

Additional difficulties arise in the blocked version that do not
appear in the single vector case. Adding single vectors to the basis
allows one to always ensures that each basis vector adds new
information. In the blocked case however, $\widehat{S}$ can become
(numerically) rank-deficient. Therefore, in
Algorithm~\ref{algo::adaptive_rs_inc},
Line~\ref{line:break_when_rank_deficient}, we look at the diagonal
elements of $\overline{R}$ in order to determine if $\widehat{S}$ is
rank-deficient. If the diagonal elements of $\overline{R}$ are less
than a specified tolerance (relative or absolute), then we have
complete knowledge of the range space (up to the specified tolerance)
and we can compress the HSS node.

%Furthermore, it is possible that the spectral decay of the off-diagonal
%blocks may not match the spectral decay of the random samples
%(\textbf{Give an example which clearly shows this}).

A similar approach is used in~\cite{Liu2016PRMatfree} with an absolute
tolerance, based on Lemma~\ref{lem:RandomLem}.
%  \todo{\chgcomment{Can we find more articles which explicitly use
% (or reference) this stopping criterion?}}
We have been unable to find an explicit
reference to a \emph{relative} stopping criterion.
From~\cite{Xia2016stability, Xia2014StableToeplitz} we know that using
a relative tolerance for the compression of off-diagonal blocks leads
to a relative error in the Frobenius norm. Since this is frequently
desired, our Algorithm~\ref{algo::adaptive_rs_inc} uses both absolute
and relative stopping criteria explicitly. Relative tolerances are
especially useful if the magnitude of different matrix sub-blocks
differ significantly. By continuing the example above, we will compute
the interpolative decomposition at a node if either condition is
satisfied:
\begin{equation}
  \label{eq:stopping_criteria}
  \frac{\min_{i}\abs{\overline{R}_{ii}}}{\rho} < \eps_{r} \, , \quad
  \min_{i} \abs{\overline{R}_{ii}} < \eps_{a} \, , \quad
  \frac{||\widehat{S}||_{F}}{||S||_{F}} < \eps_{r} \, , \quad
  \frac{1}{\sqrt{d}} ||\widehat{S}||_{F} < \eps_{a} \, .
\end{equation}
Here, $\rho$ characterizes the largest value in the $\overline{R}$
matrices. A few possible choices are $|(\overline{R_{1}})_{11}|$,
$\max_{j} |(\overline{R_{1}})_{jj}|$, or
$\max_{k,j} |(\overline{R_{k}})_{jj}|$.  As soon as one of the above
stopping criteria is satisfied, a rank-revealing factorization is
applied to all of the computed random samples.

Ideally, we would like to bound our errors with respect to $A$:
$\|(I - QQ^*)A\|_{F}/\|A\|_F$ and $\|(I - QQ^*)A\|_F$, where $A$ is a
Hankel block of the original input matrix. However, since $A$ might
not be readily available (or expensive to compute), we instead use the
random samples $S$. In the next section we establish a stochastic
$F$-norm relationship between $\|A\|_F$ and $\|S\|_F$ and, in
Appendix~\ref{sec:appendix}, show that this estimate is accurate to
high probability.

%\chgcomment{I think I should talk about $\rho$ and $R$ values,
%  possibly in this section, although this was mentioned briefly.
%  Also, I think we should talk about some of the potential
%  difficulties that arise from using non-orthonormal bases in our HSS
%  compression; see my report from September 11 (Section 6, starting at
%  page 39).}

\subsection{Mathematical Theory\label{ssec:math_theory}}
We now present the probability theory for the stopping criterion
described by Eq.~\eqref{eq:stopping_criteria}, and used in
Algorithm~\ref{algo::adaptive_rs_inc}. Let $A\in\R^{m\times n}$ and
$x\in\R^{n}$ with $x_{i}\sim\mathcal{N}(0,1)$. Let
\begin{equation}
  A = U\Sigma V^{*}
  = \begin{bmatrix} U_{1} & U_{2} \end{bmatrix}
  \begin{bmatrix} \Sigma_{r} & 0 \\ 0 & 0 \end{bmatrix}
  \begin{bmatrix} V_{1}^{*} \\ V_{2}^{*} \end{bmatrix}
\end{equation}
be the singular value decomposition (SVD) of $A$, and let
$\xi = V^{*}x$. Since $x$ is a Gaussian random vector, so is $\xi$.
By the rotational invariance of $\norm{\cdot}_{2}$, it follows that
\begin{equation}
  \norm{Ax}_{2}^{2} = \norm{\Sigma \xi}_{2}^{2}
  = \sigma_{1}^{2}\xi_{1}^{2} + \cdots + \sigma_{r}^{2}\xi_{r}^{2}\, ,
  \label{eq:prob_def}
\end{equation}
with $\sigma_{1} \ge \cdots \ge \sigma_{r} > 0$ the positive singular
values. Hence,
\begin{equation}
  \E \parens{ \norm{Ax}_{2}^{2} } = \sigma_{1}^{2} + \cdots + \sigma_{r}^{2}
  = \norm{A}_{F}^{2} \, .
\end{equation}
In order to facilitate analysis, we define
\begin{equation}
  X \sim \sigma_{1}^{2}\xi_{1}^{2} + \cdots + \sigma_{r}^{2}\xi_{r}^{2} \, ,
\end{equation}
where $\xi_{i}\sim\mathcal{N}(0,1)$, and thus
$\E \parens{X} = \norm{A}_{F}^{2}$. From~\eqref{eq:prob_def}, we see
that $\norm{Ax}_{2}^{2}$ and $X$ have the same probability
distribution, so we focus on understanding $X$. Consider
\begin{equation}
  \label{eq:averaged_RV}
  \overline{X}_{d} \sim \frac{1}{d}\parens{X_{1} + \cdots + X_{d}},
\end{equation}
%\xslcomment{Need to state that $\E (X) = \norm{A}_{F}^2$.}
where $X_{i}$ are independent realizations of $X$.  It is easy to see
that $\mathbb{E}\parens{\overline{X}_{d}} = \norm{A}_{F}^{2}$.  Using
Chernoff's inequality~\cite{prob_book}, we prove the following theorem, stating that
the probability tails of $X$ and $\overline{X}_{d}$ decay
exponentially away from $\norm{A}_F^2$.
%\xslcomment{Cite the 2 formula here: (A.9) and (A.10)?}

\begin{theorem}[Probabilistic Error Bounds]\label{thm::prob_err_bound}
  Given $\overline{X}_{d}$ as defined in Eq.~\eqref{eq:averaged_RV},
  the following bounds on the tail probabilities hold:
  \begin{align}
    \P\brackets{\overline{X}_{d}\ge \norm{A}_{F}^{2}\tau}
    &\le \exp\parens{-\frac{d\tau}{2}} \norm{A}_{F}^{dr}
      \prod_{k=1}^{r}\parens{A_{k}'}^{-d}
      \qquad \tau>1 \nonumber\\
    \P\brackets{\overline{X}_{d}\le \norm{A}_{F}^{2}\tau}
    &\le \exp\parens{-\frac{d\tau}{2}} \norm{A}_{F}^{dr}
      \prod_{k=1}^{r}\parens{A_{k}''}^{-d}
      \qquad \tau\in[0,1).
      %\vspace{-.2in}
  \end{align}
  Here, $\norm{A}_{F}^{2} = \sigma_{1}^{2} + \cdots + \sigma_{r}^{2}$,
  $\parens{A_{k}'}^{2} = \norm{A}_{F}^{2} - \sigma_{k}^{2}$, and
  $\parens{A_{k}''}^{2} = \norm{A}_{F}^{2} + \sigma_{k}^{2}$.  We know
  that $\E\parens{\overline{X}_{d}} = \norm{A}_{F}^{2}$, so $\tau$
  controls multiplicative deviation above or below the expectation
  value.
\end{theorem}
% \todo{I want to put the proof here, taken from Appendix A.1. The
% rest, A.2 and A.3 can stay in Appendix.}

% \xslcomment{This $\tau$ is unrelated to the truncation tolerance??}
% \chgcomment{Yes, $\tau$ is unrelated to the truncation tolerance.
% We can use another symbol if necessary.}

The proof is relegated to Appendix~\ref{sec:appendix}.  From
Theorem~\ref{thm::prob_err_bound}, we see that if
$R \in \R^{n \times d}$ with $R_{jk} \sim \mathcal{N}(0,1)$, it is
clear that
\begin{equation}
  \frac{1}{\sqrt{d}}\sqrt{\E\parens{\norm{AR}_{F}^{2}}} = \norm{A}_{F} \, ,
\end{equation}
and particular realizations will, with high probability, be close to
the expected value. Hence, Theorem~\ref{thm::prob_err_bound} shows
that the matrix Frobenius norm can be accurately predicted using just
(Gaussian) random samples of the range. In particular, we can
approximate the difference between our approximation and the actual
matrix sub-block, allowing us to compute both the absolute and
relative error in contrast to Lemma~\ref{lem:RandomLem}.
Future work investigating a random variable whose expectation value is
$\norm{A}_{2}$ (or some power) would be beneficial.  At this point, we
settle for using the Frobenius norm because we can accurately
approximate it.

\begin{comment}
\xslcomment{One related Slac conversion:
xiaoye:
Another question -- why is your proposed F-norm termination criterion better than the one in HMT? (see error bound (4.3) and Algorithm 4.2)
chgorman:
In HMT, the focus there is on *bounding* the 2-norm of a matrix. My estimates attempt to compute the expectation value of the 2-norm from many samples. In HMT, there bounds are guaranteed, but with a probability of failure that can be made negligibly small with the potential for (large) error overestimates. Right now, mine only says that what our error estimates will hover around. More theoretical work will be required to give a confidence interval for individual realizations. I think this can be done, but right now I am not sure what will be required.
The proposed F-norm termination criterion is potentially better because we should more closely match the actual error than overestimate it. The downside is that we do not bound the matrix 2-norm but rather the matrix F-norm.
}
\end{comment}

\subsection{Flop Counts: {\dbl} vs. {\inc}}
First, note that the full QR factorization $Q = \text{QR}(A)$ for
$A\in\mathbb{R}^{m\times n}$ (with $m \gg n$) performs $2mn^{2}$ floating point
operations (flops).  Assuming the numerical rank of $A$ is $r$, then
the rank-revealing QR factorization $Q = \text{RRQR}(A)$ requires
$2mnr$ flops. %% \xslcomment{double check??}
Given $S\in\mathbb{R}^{m\times r}$ and $X\in\mathbb{R}^{m\times d}$,
the orthogonalization $S \leftarrow (I - X X^{*})S$ requires $4mdr$
flops (ignoring a lower order $mr$ term), and twice that for the
iterated ($2\times$) block Gram-Schmidt step.
% The reorthogonalization is the most expensive part of the algorithm,
% so we will focus on this next.

In the {\dbl} strategy, Algorithm~\ref{algo::adaptive_rs_dbl}, at step
$k$, $2^{k-1} d_0 + p$ random vectors have already been sampled, where
$p$ is a small oversampling parameter (e.g. $p = 10$).  Then,
$2^{k-1}d_0$ new random samples are added to the sample matrix $S_k$,
leading to $2^{k} d_0 + p$ columns for $S_k$.
% and compute $S_k$ with $2^{k-1} d_0$ columns.
Except for the final step, the operation
$Q_{k} \leftarrow \text{RRQR}(S_{k})$ costs $\sim 2 m 4^{k} d_0^{2}$
flops. In the final step, the RRQR terminates early when the tolerance
is met at rank $r$, with a cost of $\sim 2 m r^2$ flops. 
%\todo{count final step separately}
The total number of steps needed to reach the
final rank $r$ is $N\sim \log(r / d_0)$.  Summing the costs of $N$ steps:
\(2 m d_0^2 \sum_{k=1}^{N}4^k \), we obtain the total flop count of
$\sim \frac{8}{3} mr^{2}$.
%\xslcomment{Porbably can comment on the best / worst cases.}
\ignore{ \xslcomment{I do not understand the following: ``In the worst
    case, we see that the worst case flop count is $\sim 8mr^{2}$.
    Here, we see that we do four times more work in the work case
    asymptotically. The cost of RRQR is $2mr^{2}$ flops, wtih
    $8mr^{2}$ flops in the worst case.  Asymptotically, the cost is
    $4mr^{2}$ flops in the best case and $16mr^{2}$ flops in the worst
    case.''  } }

The {\inc} strategy, Algorithm~\ref{algo::adaptive_rs_inc}, starts
with $d_{0}$ random samples and adds $\Delta d$ new random vectors at
each step. At step $k$ we have the sample matrix
$S = [ S_{1} \,\, \cdots \,\, S_{k} ]$, and the orthogonal matrix
$Q = [Q_{1} \,\, \cdots \,\, Q_{k-1}]$.  We first compute the
orthogonal projection $\hat{S}_k \gets (I - Q Q^{*})^2S_k$, which
costs $8m (k-1) \Delta d^2$ flops, followed by
$Q_{k} \gets \text{QR}(\hat{S}_{k})$, which costs $2 m \Delta d^2$,
and then append $Q_k$ to $Q$.  The total number of steps needed to
reach the final rank $r$ is
$N \sim (r - d_{0}) / \Delta d \sim r/\Delta d$.  Summing the cost of
$N$ steps: \( 8 m {\Delta d}^2 \sum_{k=1}^N (k-1) \), gives the
overall cost $\sim 4 mr^{2}$.
%The worst-case flop count is asymptotically the same as {\dbl}.
The cost of the final RRQR is an additional $2mr^{2}$ flops.
% \xslcomment{I got $2mr \Delta d$.}
% Asymptotically, the total cost is $4mr^{2}$ flops.

%\input{tex/flop_count.tex}

From the above analysis, we see that the {\inc} scheme requires more
flops. However, {\dbl} involves a larger amount of data movement due
to the column pivoting needed in each step of RRQR. This manifests
itself in the communication cost of the parallel algorithm, which will
be analyzed in Section~\ref{sec:comm_cost}.

%% file: tex/tables/defs.tex
%%%%%%%%%%%%%%%%%%%%%%%%%%%%%%%%%%%%%%%%%%%%%%%%%%%%%%%%%%%%%%%%%%%%%%%%
%%% Definition Table

\begin{table}
  \begin{center}
    \begin{tabular}[c]{r|l}
      \hline
      \texttt{randn($m, n$)}    & an $m \times n$ matrix with iid $\mathcal{N}(0, 1)$ elements \\
      \texttt{rows($A$)/cols($A$)} & number of rows/columns in matrix $A$ \\
      %\texttt{local($\tau$)}    & \texttt{true} if this process owns or works on node $\tau$, \texttt{false} otherwise \\
      \texttt{isleaf($\tau$)}   & \texttt{true} if $\tau$ is a leaf node, \texttt{false} otherwise \\
      \texttt{isroot($\tau$)}   & \texttt{true} if $\tau$ is a root node, \texttt{false} otherwise \\
      \texttt{children($\tau$)} & a list with the children of node $\tau$, always zero or two \\
      \texttt{levels($\mathcal{T}$)} & number of levels in tree $\mathcal{T}$ \\
      % \texttt{level($\mathcal{T}$, $\ell$)} & all nodes of tree $\mathcal{T}$ at level $\ell$, level $1$ is the root \\
      \texttt{level($\tau$)} & level of node $\tau$, starting from $0$ at the root \\
      \texttt{$\{Q, r\} \gets$ RRQR($S, \varepsilon_r, \varepsilon_a$)} & rank-revealing QR, orthonormal $Q$ and rank $r$\\
      % \texttt{allgatherv($\{\mathcal{J}_1, \dots, \mathcal{J}_P\}$)} & every process $p \in 1 : P$ broadcasts $\mathcal{J}_p$ to all processes $\in 1:P$ \\
      %\texttt{send(d, $\mathcal{I}$)} & send $\mathcal{I}$ to processor \texttt{d} \\
      %\texttt{receive()}        & receive a message \\
      %\texttt{nrprocs(comm)}    & number of processes in communicator \texttt{comm} \\
      %\texttt{rank(comm)}       & rank of this process in communicator \texttt{comm} \\
      \hline
    \end{tabular}
  \end{center}
  \caption{List of helper functions.}
  \label{tab::addfunctions}
\end{table}

%% file: tex/algo/hss_comp.tex
%%%%%%%%%%%%%%%%%%%%%%%%%%%%%%%%%%%%%%%%%%%%%%%%%%%%%%%%%%%%%%%%%%%%%%%%
%%% HSS Compression Algorithm

\begin{algorithm2e}
  \DontPrintSemicolon
  \SetAlgoLined
  \SetKwProg{Fn}{function}{}{}
  \SetKwFunction{compress}{HSSCompress}
  \SetKwFunction{compressnode}{CompressNode}
  \SetKwFunction{root}{root}
  \SetKwFunction{rows}{rows}
  \SetKwFunction{cols}{cols}
  \SetKwFunction{randn}{randn}
  %\SetKwFunction{extract}{ExtractElements}
  \SetKwFunction{localsamples}{ComputeLocalSamples}
  \SetKwFunction{reducelocalsamples}{ReduceLocalSamples}
  \SetKwFunction{isleaf}{isleaf}
  \SetKwFunction{isroot}{isroot}
  \SetKwFunction{level}{level}
  \SetKwFunction{child}{children}
  \SetKwFunction{local}{local}
  \SetKwFunction{ID}{ID}

  \Fn{$H =$ \compress{$A$, $\mathcal{T}$}}{
    $R \gets$ \randn{\rows{$A$}, $d_0$} \\
    \compressnode{$A$, $R$, $AR$, $A^*R$, \root{$\mathcal{T}$}} \\
    \Return $\mathcal{T}$
  }
  \BlankLine
  \Fn{\compressnode{$A$, $R$, $S^r$, $S^c$, $\tau$}} {
    \uIf{\isleaf{$\tau$}}{
      $\tau.D \gets A(\tau.I, \tau.I)$ \\
    }
    \Else{
      $\nu_1, \nu_2 \gets$ \child{$\tau$} \\
      \compressnode{$A$, $R$, $S^r$, $S^c$, $\nu_1$} \tcp*{recursive compression ..}
      \compressnode{$A$, $R$, $S^r$, $S^c$, $\nu_2$} \tcp*{of child nodes}
      $\tau.B_{12} \gets A(\nu_1.I^r, \nu_2.I^c)$ \\
      $\tau.B_{21} \gets A(\nu_2.I^r, \nu_1.I^c)$ \\
    }
    \lIf{\isroot{$\tau$}}{\Return}
    \localsamples{$R$, $S^r$, $S^c$, $\tau$, $1:\cols{R}$} \\
    $\{ (\tau.U)^*,\,\, \tau.J^r \} \gets$ \ID{$(\tau.S^r)^*$,
      $\varepsilon_r / $\level{$\tau$}, $\varepsilon_a / $\level{$\tau$}} \\
    $\{ (\tau.V)^*,\,\, \tau.J^c \} \gets$ \ID{$(\tau.S^c)^*$,
      $\varepsilon_r / $\level{$\tau$}, $\varepsilon_a / $\level{$\tau$}} \\
    \reducelocalsamples{$R$, $\tau$, $1:\cols{R}$}
  }
  \caption{Non-adaptive HSS compression of
    $A \in \mathbb{R}^{N \times N}$ using cluster tree $\mathcal{T}$
    with relative and absolute tolerances $\varepsilon_r$ and
    $\varepsilon_a$ respectively. }
  \label{algo:HSScompress}
\end{algorithm2e}

%% file: tex/algo/comp_local_reduce_2.tex
%%%%%%%%%%%%%%%%%%%%%%%%%%%%%%%%%%%%%%%%%%%%%%%%%%%%%%%%%%%%%%%%%%%%%%%%
%%% Compute Local Samples and Reduce Samples 2

\begin{algorithm2e}  %[!ht]
%\begin{algorithm}[p]
  \DontPrintSemicolon
  \SetAlgoLined
  \SetKwProg{Fn}{function}{}{}
  \SetKwFunction{localsamples}{ComputeLocalSamples}
  \SetKwFunction{isleaf}{isleaf}
  \SetKwFunction{cols}{columns}
  \SetKwFunction{child}{children}
  \SetKwFunction{reducelocalsamples}{ReduceLocalSamples}

  \Fn{\localsamples{$R$, $S^r$, $S^c$, $\tau$, $i$}} {
    % $i \gets a : b$ \tcp*[f]{Index set of only the last $\Delta d$ columns} \\
    \uIf{\isleaf{$\tau$}}{
      $\tau.S^r(\,:\,,i) \gets S^r(\tau.I,i) - \tau.D \,\, R(\tau.I,i)$ \\
      $\tau.S^c(\,:\,,i) \gets S^c(\tau.I,i) - (\tau.D)^* \,\, R(\tau.I,i)$ \\
    }\Else{
      $\nu_1, \nu_2 \gets$ \child{$\tau$} \\
      $\tau.S^r(\,:\,,i) \gets \begin{bmatrix} \nu_1.S^r(\nu_1.J^r,i) - \tau.B_{12} \,\, \nu_2.R^r \\
        \nu_2.S^r(\nu_2.J^r,i) - \tau.B_{21} \,\, \nu_1.R^r \end{bmatrix}$ \\
      $\tau.S^c(\,:\,,i) \gets \begin{bmatrix} \nu_1.S^c(\nu_1.J^c,i) - (\tau.B_{21})^* \,\, \nu_2.R^c \\
        \nu_2.S^c(\nu_c.J^c,i) - (\tau.B_{12})^* \,\, \nu_1.R^c \end{bmatrix}$ \\
    }
  }
  \BlankLine
  \Fn{\reducelocalsamples{$R$, $\tau$, $i$}}{
    % $\tau.S^r(\,:\,,i) \gets \tau.S^r(\tau.J^r, i)$ \\
    % $\tau.S^c(\,:\,,i) \gets \tau.S^c(\tau.J^c, i)$ \\
    \uIf{\isleaf{$\tau$}}{
      $\tau.R^r \gets (\tau.V)^* \,\, R(\tau.I,i)$ ; \quad
      $\tau.I^c \gets \tau.I(\tau.J^c)$ \\
      $\tau.R^c \gets (\tau.U)^* \,\, R(\tau.I,i)$ ; \quad
      $\tau.I^r \gets \tau.I(\tau.J^r)$
    }\Else{
      $\nu_1, \nu_2 \gets$ \child{$\tau$} \\
      $\tau.R^r \gets (\tau.V)^* \begin{bmatrix} \nu_1.R^r \\ \nu_2.R^r \end{bmatrix}$ ; \quad
      $\tau.I^r \gets \begin{bmatrix} \nu_1.I^r & \nu_2.I^r \end{bmatrix}(\tau.J^r)$ \\
      $\tau.R^c \gets (\tau.U)^* \begin{bmatrix} \nu_1.R^c \\ \nu_2.R^c \end{bmatrix}$ ; \quad
      $\tau.I^c \gets \begin{bmatrix} \nu_1.I^c & \nu_2.I^c \end{bmatrix}(\tau.J^c)$
    }
  }
  \caption{Compute local samples and reduce local samples based on
    rows selected by the interpolative decomposition.}
  \label{algo::computelocalsamples}
\end{algorithm2e}

%% file: tex/algo/adaptive_comp_node.tex
%%%%%%%%%%%%%%%%%%%%%%%%%%%%%%%%%%%%%%%%%%%%%%%%%%%%%%%%%%%%%%%%%%%%%%%%
%%% Adaptive Node Compression

\begin{algorithm2e}
  \DontPrintSemicolon
  \SetAlgoLined
  \SetKwProg{Fn}{function}{}{}
  \SetKwFunction{compressadaptive}{HSSCompressAdaptive}
  \SetKwFunction{compressnodeadaptive}{CompressNodeAdaptive}
  \SetKwFunction{localsamples}{ComputeLocalSamples}
  \SetKwFunction{reducelocalsamples}{ReduceLocalSamples}
  \SetKwFunction{isleaf}{isleaf}
  \SetKwFunction{isroot}{isroot}
  \SetKwFunction{child}{children}
  \SetKwFunction{ID}{ID}
  \SetKwFunction{local}{local}
  \SetKwBlock{Try}{try}{}
  \SetKwBlock{Fail}{catch}{}
  \SetKwData{UNTOUCHED}{\texttt{UNTOUCHED}}
  \SetKwData{COMPRESSED}{\texttt{COMPRESSED}}
  \SetKwData{PARTIALLYCOMPRESSED}{\texttt{PARTIALLY\_COMPRESSED}}
  \SetKwData{state}{state}

  \Fn{$H =$ \compressadaptive{$A$, $\mathcal{T}$}}{
    $d \gets d_0$; \quad $\Delta d \gets 0$; \quad
    $R \gets $ \randn{$N, d$}; \quad
    $S^r \gets AR$; \quad
    $S^c \gets A^*R$ \\
    \lForEach{$\tau \in \mathcal{T}$}{$\tau.\state \gets \UNTOUCHED$}
    \While{\root{$\mathcal{T}$}$.\state \neq \COMPRESSED$ {\bf and} $d < d_{\text{max}}$}{
      \compressnodeadaptive{$A$, $R$, $S^r$, $S^c$, \root{$\mathcal{T}$}, $\varepsilon$, $d$, $\Delta d$}\\
      %\If{\root{$\mathcal{T}$}$.\state \neq \COMPRESSED$ {\bf and} $d < d_{\text{max}}$}{
      $\Delta d \gets d \,\, \text{or} \,\, c^{st}$ \tcp*{double or increment?}
      $\bar{R} \gets$ \randn{$N$, $\Delta d$}; \quad
      $R \gets \begin{bmatrix} R & \bar{R} \end{bmatrix}$ \\
      $S^r \gets \begin{bmatrix} S^r & A \bar{R} \end{bmatrix}$; \quad
      $S^c \gets \begin{bmatrix} S^r & A^* \bar{R} \end{bmatrix}$ \\
      $d \gets d + \Delta d$ \\
      %}
    }
    \Return $\mathcal{T}$
  }
  \BlankLine
  \Fn{\compressnodeadaptive{$A$, $R$, $S^r$, $S^c$, $\tau$, $d$, $\Delta d$}} {
    \uIf{\isleaf{$\tau$}}{
      \lIf{$\tau.\state = \UNTOUCHED$}{
        $\tau.D \gets A(\tau.I, \tau.I)$
      }
    }
    \Else{
      $\nu_1, \nu_2 \gets$ \child{$\tau$} \\
      \compressnodeadaptive{$A$, $R$, $S^r$, $S^c$, $\nu_1$, $d$, $\Delta d$} \\
      \compressnodeadaptive{$A$, $R$, $S^r$, $S^c$, $\nu_2$, $d$, $\Delta d$} \\
      \lIf{$\nu_1.\state \neq \COMPRESSED$ {\bf or} $\nu_2.\state \neq \COMPRESSED$}{\Return}
      \uIf{$\tau.\state = \UNTOUCHED$}{
        $\tau.B_{12} \gets A(\nu_1.I^r, \nu_2.I^c)$; \quad
        $\tau.B_{21} \gets A(\nu_2.I^r, \nu_1.I^c)$
      }
    }
    \lIf{\isroot{$\tau$}}{
      $\tau.\state \gets \COMPRESSED$; \quad
      \Return
    }
    \uIf{$\tau.\state = \UNTOUCHED$}{
      \localsamples{$R$, $S^r$, $S^c$, $\tau$, $1:d+\Delta d$} %\tcp*[f]{Compute all $d$ columns of $\tau.S^{r/c}$ \& $\tau.R^{r/c}$} \pgcomment{arguments!} \\
    }\lElse{
      \localsamples{$R$, $S^r$, $S^c$, $\tau$, $d+1:d+\Delta d$} %\tcp*[f]{Only last $\Delta d$ columns}
    }
    \uIf{$\tau.\state \neq \COMPRESSED$}{
      \Try{
        $\{ (\tau.U)^*,\,\, \tau.J^r \} \gets$ \ID{$(\tau.S^r)^*$,
          $\varepsilon_r / $\level{$\tau$}, $\varepsilon_a / $\level{$\tau$}} \\
        $\{ (\tau.V)^*,\,\, \tau.J^c \} \gets$ \ID{$(\tau.S^c)^*$,
          $\varepsilon_r / $\level{$\tau$}, $\varepsilon_a / $\level{$\tau$}} \\
        $\tau.\state \gets \COMPRESSED$ \\
        \reducelocalsamples{$R$, $\tau$, $1:d+\Delta d$} %\tcp*[f]{Update all $d$ columns of $\tau.S^{r/c}$ \& $\tau.R^{r/c}$}
      }
      \Fail(\tcp*[f]{RRQR/ID failed to reach tolerance $\varepsilon_r$ or $\varepsilon_a$}){
        $\tau.\state \gets \PARTIALLYCOMPRESSED$ \\
        \Return
      }
    }\lElse{
      \reducelocalsamples{$R$, $\tau$, $d$, $d+1:d+\Delta d$} %\tcp*[f]{Update last $\Delta d$ columns of $\tau.S^{r/c}$ \& $\tau.R^{r/c}$}
    }
  }
  \caption{Adaptive HSS compression of
    $A \in \mathbb{R}^{N \times N}$ using cluster tree $\mathcal{T}$
    with relative and absolute tolerances $\varepsilon_r$ and
    $\varepsilon_a$ respectively. }
  \label{algo::HSScompressNodeAdaptive}
\end{algorithm2e}

%% file: tex/algo/adaptive-rs-dbl.tex
%%%%%%%%%%%%%%%%%%%%%%%%%%%%%%%%%%%%%%%%%%%%%%%%%%%%%%%%%%%%%%%%%%%%%%%%
%%% Doubling strategy

\begin{algorithm2e}
  \DontPrintSemicolon
  \SetAlgoLined
  \SetKwProg{Fn}{function}{}{}
  \SetKwFunction{randn}{randn}
  \SetKwFunction{ID}{ID}
  \SetKwFunction{rsdbl}{RS-\dbl}
  \SetKwFunction{RRQR}{RRQR\_HMT}

  \Fn{$Q =$ \rsdbl{$A$, $d_0$, $p$, $\eps_r$, $\eps_a$}}{
    $k\gets 1$; \quad $m\gets \rows{A}$; \quad $r \gets \infty$ \\
    $R_1 \gets $ \randn{$m$, $d_0 + p$} \\
    $S_1 \gets A R_1 $ \\
    \While{$(r > 2^{k-1} d_0)$} {
      $\{Q, r\} \gets$ \RRQR{$S_k$, $\eps_r$, $\eps_a$} \\
      $\Delta d \gets 2^{k-1} d_0 + p$  \tcp*{double sample size}
      $R_{k+1} \gets$ \randn{$m$, $\Delta d$} \\
      $S_{k+1} \gets [S_k \  A R_{k+1} ]$ \\
      $k \gets k + 1$ \\
    }
    \Return $Q$
  }
  \caption{Adaptive computation of $Q$, an approximate basis for the
    range of the Hankel block $A$, using the {\dbl} strategy with an
    oversampling parameter $p$.}
  \label{algo::adaptive_rs_dbl}
\end{algorithm2e}

%% file: tex/algo/adaptive-rs-inc.tex
%%%%%%%%%%%%%%%%%%%%%%%%%%%%%%%%%%%%%%%%%%%%%%%%%%%%%%%%%%%%%%%%%%%%%%%%
%%% Increment strategy

\begin{algorithm2e}
  \DontPrintSemicolon
  \SetAlgoLined
  \SetKwProg{Fn}{function}{}{}
  \SetKwFunction{randn}{randn}
  \SetKwFunction{ID}{ID}
  \SetKwFunction{rsinc}{RS-\inc}
  \SetKwFunction{RRQR}{RRQR}
  \SetKwFunction{QR}{QR}
  \SetKwFunction{break}{break}
  \SetKwFunction{min}{min}
  \SetKwFunction{diag}{diag}

  \Fn{$Q =$ \rsinc{$A$, $d_0$, $\Delta d$, $\eps_r$, $\eps_a$}}{
    $k\gets 1$; \quad $m\gets \rows{A}$; \quad $Q \gets [ \ ]$ \\
    $R_1 \gets$ \randn{$m$, $d_0$} \\
    \While{true}{
      $S_k \gets A R_k$ \tcp*{new samples}
      %\pgcomment{do this twice? added $^2$} \\
      $\hat{S}_k \gets (I - QQ^*)^2 S_k$ \tcp*{iterated block Gram-Schmidt} \label{line:block_Gram_Schmidt}
      % \tcp*{ortho.~against previous $Q$}
      \lIf{$(||\hat{S}_k||_{F} / ||S_k||_{F} < \eps_{r}$
        or $||\hat{S}_k||_{F} / \sqrt{\Delta d} < \eps_{a})$} {
        \break
      }
      $\{Q_k, \overline{R}_{k}\} \gets$ \QR{$\hat{S}_k$} \label{line:inc_QR} \\
      \lIf{\min{\diag{$|\overline{R}_{k}|$}}$ < \eps_a$ or
        \min{\diag{$|\overline{R}_{k}|$}}$ < \eps_r |(\overline{R}_{1})_{11}|$} {
        \break \label{line:break_when_rank_deficient}
      }
      %\chgcomment{Include because present in algorithm?} \\
      $Q \gets [ Q \ Q_k ]$ \\
      $k \gets k + 1$ \\
      $R_k \gets$ \randn{$m$, $\Delta d$} \\
    }
    %\pgcomment{RRQR of Q or of S??} \\
    %$\{Q, r\} \gets$ \RRQR{$Q$, $\eps_r$, $\eps_a$} \\
    $\{Q, r\} \gets$ \RRQR{$[S_1 \ \dots \ S_k]$, $\eps_r$, $\eps_a$} \label{line::incrementing_RRQR}\\
    %\xslcomment{Sufficient to do only $RRQR(Q_k)$?} \\
    \Return $Q$
  }
  \caption{Adaptive computation of $Q$, an approximate basis for the
    range of the Hankel block $A$, using the {\inc} strategy.}
  \label{algo::adaptive_rs_inc}
\end{algorithm2e}

%% file: tex/parallel_algorithm.tex
\section{Parallel Algorithm\label{sec:parallel}}
The parallel algorithm uses the same parallelization framework as
described in~\cite[Section 3]{FHRdistributedHSS}. The data
partitioning and layout is based on the HSS tree, following a
proportional mapping of subtrees to subsets of processes in a top-down
traversal. The HSS tree can be specified by the user. The tree should
be binary, but can be imbalanced and does not need to be complete. For
HSS nodes that are mapped to multiple processors, the matrices stored
at those nodes are distributed in 2D block-cyclic (ScaLAPACK style)
layout.

%\subsection{Data Partitioning and Layout}

\subsection{Partially Matrix-Free Interface}
The randomized HSS compression algorithm is so-called ``partially''
matrix-free. This means it does not require every single element of
the input matrix $A$. What is required is a routine to perform the
random sampling $S = AR$, as well as a way to extract sub-blocks from
$A$. Recall from Section~\ref{ssec:randomHSS} that at the leafs,
$D_\tau = A(I_\tau,I_\tau)$. Furthermore, due to the use of
interpolative decompositions,
$B_{\nu_1,\nu_2} = A(J_{\nu_1}^r, J_{\nu_2}^c)$, is a sub-block of
$A$.

If the input matrix $A$ is given as an explicit dense matrix, the
random sampling $S=AR$ is performed in parallel using the PBLAS
routine \texttt{PDGEMM} with a 2D block-cyclic data layout for $A$,
$S$ and $R$. In this case, the input matrix $A$ is also redistributed
-- with a single collective MPI call -- from the 2D block-cyclic
layout to a layout corresponding to the HSS tree, such that extraction
of sub-blocks for $D_\tau$ and $B_{\nu_1,\nu_2}$ does not require
communication between otherwise independent HSS subtrees.

Instead of forming an explicit dense matrix $A$, the user can also
specify multiplication and element extraction routines. The
multiplication routine computes, for a given random matrix $R$, the
random sample matrices $S^r = A R$ and $S^c = A^* R$. This is the more
interesting use case, since for certain classes of structured matrices
a fast multiplication algorithm is available; consider for instance
sparse matrices, low-rank and hierarchical matrices, combinations of
sparse and low-rank matrices or operators which can be applied using
the fast Fourier transform or similar techniques. The element
extraction routine should be able to return matrix sub-blocks
$A(I,J)$, defined by row and column index sets $I$ and $J$
respectively. Depending on how the user data is distributed, computing
matrix elements might involve communication between all processes. In
this case, the HSS compression traverses the tree level by level (from
the leafs to the root), with synchronization at each level and element
extraction for all blocks $D_\tau$ and $B_{\nu_1,nu_2}$ on the same
performed simultaneously in order to aggregate communication messages
and minimize communication latency. If no communication is required
for element extraction, then independent subtrees can be compressed
concurrently.
% \todo[inline]{say something about element extraction in this case, it
%   still works if it requires communication.}

% We require matrix-vector multiplication and element extraction
% routine.

\ignore{ %%%% ignored 
\subsection{Synchronized Element Extraction}\label{sec:sync_extract}
Synchronization of element extraction occurs when the information
required to extract matrix elements is distributed between
processors so that they must work together.
%Add synchronization to the algorithm

Algorithm~\ref{algo::HSScompressNodeSync}.
%%% TODO %%%% convert to algorithmic?? %%%%%
%%\input{tex/algo/comp_node_sync.tex}
%%% TODO %%%%%%%%%%%%%%%%%%%%%%%%%%%%%%%%%%%

TODO communicate states of nodes to all processes, communicate state
up the tree??

Algorithm~\ref{algo::exchange_indices}
%%% TODO %%%% convert to algorithmic?? %%%%%
%%\input{tex/algo/exchange_indeces.tex}
%%% TODO %%%%%%%%%%%%%%%%%%%%%%%%%%%%%%%%%%%

Synchronization for element extraction can be done in two ways: 1)
Participate\_to\_calls with send/receive to communicate I,J 2)
allgather of I,J
\\
% No need for the following (too complicated to implement)
% In the end, use a hybrid strategy: use send/receive on higher levels
% in the HSS tree, allgather below
What is communication cost?
\xslcomment{This part used to dominate runtime on large core count.
  See IPDPS17 paper.  How much is improved now?}
}

\subsection{Parallel Restart}
During factorization, nodes can be in either \texttt{UNTOUCHED}
(\texttt{U}), {\PC} (\texttt{PC}) or \texttt{COMPRESSED} (\texttt{C})
state. A node can not start compression until both it's children are
in the \texttt{C} state. If during HSS compression, a process
encounters an internal HSS node which children that are not both in
the \texttt{C} state, this process stops the HSS tree
traversal. Independent subtrees can progress the compression further
if they can compressed successfully with the current number of random
samples. This can lead to load imbalances if the HSS tree, or the
off-diagonal block ranks, are imbalanced. Eventually, all processes
synchronize to perform the random sampling in parallel. Hence there is
some overhead associated with restarting the HSS compression algorithm
to add more random samples. In addition, random sampling is more
efficient, in terms of floating point throughput, when performed with
more random vectors at once.

% Draw a picture of HSS tree to illustrate the moment when restart is
% needed.  Show the node states:

\subsection{Communication Cost in Parallel Adaptation}
\label{sec:comm_cost}
%  {\dbl} has higher communication cost.
In~\cite{FHRdistributedHSS}, we analyzed the communication cost of the entire
parallel HSS algorithm, assuming no adaptivity. In this section,
we will focus only on the cost of adaptivity, using either {\dbl} or
{\inc} strategy.

Consider the current node of the HSS tree that requires adaptation
(called ``PC'' node). Assume the final rank is $r$, the row dimension of the
sample matrix is $m$, and $P$ processes work on this node in parallel. 
We use the pair [\#messages, volume] to denote the {\em communication cost}
which counts the number of messages and the number of words transferred
for a given operation, typically along the critical path. 
A broadcast of $W$ words among $P$ processes has the cost
[$\log P$, $W \log P$]. This assumes that broadcast follows a tree-based
implementation: there are $\log P$ steps on the critical path
(any branch of the tree) and $W$ words are transferred at each step,
yielding $\log P$ messages and $w \log P$ words.

For a rectangular matrix of dimension $m\times n$, assuming $m/P\ge n$, 
the communication cost for non-pivoted QR factorization
(\texttt{PDGEQRF} in ScaLAPACK) is \par \noindent
[$2n\log P, \frac{m n}{\sqrt{P}}\log P$]~\cite{scalapackmanual,OptQRLU2012}.
For QR factorization with column pivoting, i.e., 
\texttt{PDGEQPF} in ScaLAPACK,
 additional communication is needed at each step to
compute the column norm and permutes the column with the maximum norm to
the leading position. Computing the maximum column norm needs two reductions
along row and column dimensions, costing $2 \log \sqrt{P} = \log P$ messages.
The additional communication volume is of lower order term. In total,
\texttt{PDGEQPF} has communication cost [$3n\log P, \frac{mn}{\sqrt{P}}\log P$]

In the {\dbl} strategy (Algorithm~\ref{algo::adaptive_rs_dbl}),
we need {$s = \log{\frac{r}{d_0}}$} steps of augmentations
to reach the final rank $r$. At the $k$-th step,
we perform RRQR (ID) for $S_k$ of dimension $m\times d_0 2^k$,
using \texttt{PDGEQPF}. The total communication cost of $s$ steps sums up to:
\begin{align}\label{eq:dbl_comm_cost}
\sum_{k=0}^{s} \left[ 3 d_0 2^k \log P,\; \frac{m (d_0 2^k)}{\sqrt{P}}\log P \right] =  
\left[ 3 d_0 \log P\sum_{k=0}^{s} 2^k,\; \frac{m d_0}{\sqrt{P}}\log P \sum_{k=0}^{s} 2^k\right ] \\ \nonumber
= \left[ 3 d_0 2^{s+1} \log P,\: \frac{m d_0 2^{s+1}}{\sqrt{P}} \log P \right ]
= {\left[6r \log P,\; \frac{2 m r}{\sqrt{P}}\log P \right]}.
\end{align} 

In the {\inc} strategy (Algorithm~\ref{algo::adaptive_rs_inc}),
we need {$s = \frac{r}{\Delta d}$} steps of increments
 to reach the final rank $r$. At the $k$-th step, two costly operations
are the block Gram-Schmidt and block QR 
(Lines~\ref{line:block_Gram_Schmidt} and~\ref{line:inc_QR} respectively in Algorithm~\ref{algo::adaptive_rs_inc}).

Each Gram-Schmidt orthogonalization step requires two matrix multiplications;
we use \texttt{PDGEMM} in PBLAS, which uses a pipelined SUMMA
algorithm~\cite{summa}.
The $Q$ matrix is of dimension $m\times d$, where $d = \Delta d (k-1)$.
The $S_k$ matrix is of dimension $m\times \Delta d$.
The communication cost for $Q^*\cdot S_k$ is
$[\Delta d, \frac{\Delta d (m + d)}{\sqrt{P}}]$.
The cost for another multiplication $Q\cdot (Q^*\cdot S_k)$ is
the same. Since we do Gram-Schmidt twice, the total communication
cost of $s$ steps sums up to:
\begin{equation}\label{eq:GS_comm_cost}
2 \cdot 2 \sum_{k=1}^{s} \left[ \Delta d, \frac{ \Delta d (m+d)}{\sqrt{P}}  \right] = 
4 \cdot \left[ r, \sum_{k=1}^{s} \frac{m\Delta d + {\Delta d}^2 \cdot k}{\sqrt{P}} \right] = 
\left[ 4 r, 4 \frac{m r + r^2/2}{\sqrt{P}} \right].
\end{equation}

Next, we perform QR for $\hat{S}_k$ of dimension $m\times \Delta d$,
using \texttt{PDGEQRF}. The total communication cost of $s$ steps sums up to:

\begin{equation} \label{eq:QR_comm_cost}
\sum_{k=1}^{s} \left[ 2\: \Delta d \log P, \; \frac{m\: \Delta d}{\sqrt{P}}\log P \right]
 = \left[2 r \log P, \; \frac{m r}{\sqrt{P}} \log P \right].
\end{equation}

Comparing Eqs.~\eqref{eq:GS_comm_cost} and~\eqref{eq:QR_comm_cost}, we see that
the communication in Gram-Schmidt is a lower order term compared to the QR
factorizations, therefore we ignore it.

In a final step (Line~\ref{line::incrementing_RRQR} in Algorithm~\ref{algo::adaptive_rs_inc}) 
we perform a large RRQR, with communication cost
[$3 r\log P, \; \frac{m r}{\sqrt{P}} \log P $].
This is added into Eq.~\eqref{eq:QR_comm_cost}. 
The leading cost of the {\inc} strategy is 
$[5 r\log P, \; \frac{2 m r}{\sqrt{P}} \log P ]$.
Compared to the {\dbl} strategy, the {\inc} strategy requires fewer
messages, and has a similar communication volume.

%% The cost of ``Compute Samples'', ``Reduce Samples'' are smaller, 
%% hence, we ignore them here.

%% file: tex/experiments.tex
%%%%%%%%%%%%%%%%%%%%%%%%%%%%%%%%%%%%%%%%%%%%%%%%%%%%%%%%%%%%%%%%%%%%%%%%
%%% Numerical Experiments

\section{Numerical Experiments}
\label{sec:experiments}

\subsection{Experiments Setup}

Numerical experiments were conducted using the Cori supercomputer at NERSC. Each node has two sockets, each socket is populated with a 16-core Intel\textcopyright~Xeon\texttrademark~Processor E5-2698 v3 (``Haswell'') at 2.3 GHz and 128 GB of RAM memory. We used STRUMPACK v2.2 linked with Intel MKL and compiled with the Intel compilers version 17.0.1.132 in Release mode.

\subsection{Test Problems}

The first set of numerical experiments depicts six dense linear systems from a variety of applications in double precision (dp) and single-precision (sp). A complete description of the problems under consideration can be found on \cite{FHRdistributedHSS}, with the exception of the first experiment at Table \ref{table:sixdrivers}. The first experiment is the parameterized example $\alpha I + \beta UDV^{*}$. 
$U$ and $V$ are orthogonal matrices of rank $r$ and are distributed
in order to allow for fast element extraction and scalable tests,
and $D$ is either the identity matrix or a diagonal matrix with decaying entries, giving us a matrix that has off-diagonal blocks with either constant or decaying singular values.
In this particular case we set $D_{k,k} = 2^{-53(k-1)/r}$, and $r=500$.
The largest experiments ($N=500,000$) used $p=1,024$ cores, whereas the rest of the experiments used $p=64$ cores. Column three, HSS$_\epsilon$, denotes the relative error of the HSS approximation, whereas the absolute error is kept constant at $1E-08$. Each experiment is performed at three increasingly tighter approximation tolerances, while we report on the three main stages involved in the solution of the linear system: Approximation --compression-- of the dense matrix, factorization, and solve. The metrics of interest at each stage are memory
consumption, flops, and max wall-clock time from all MPI ranks.

Numerical experiments (Table \ref{table:sixdrivers}) show that compression takes the most flops and time, followed by factorization and solve, which are much cheaper in comparison. Nonetheless, as we increase the accuracy of the approximation, we notice a moderate increase in the wall-clock time in all three stages.

\input{tables/table6drivers.tex}

\subsection{Performance Breakdown}

The second set of numerical experiments illustrates the new {\inc}
adaptive technique in comparison with the traditional {\dbl} adaptive
technique that relies on RRQR with column pivoting.
Table~\ref{table:PerformanceBreakdown} shows the detailed breakdown
of the flops and time in different stages of the two algorithms
for the BEM Spehere linear system with $d_0 = 128$  and $\Delta d = 256$,
and using $p = 1,024$ cores.

The first observation is that in both {\inc} and {\dbl}, nearly
90\% of the flops are in the initial sampling step. On the other hand,
since the sampling step involves highly efficient matrix-matrix multiplication,
the percentage time spent at this step is under 50\%.

The second most costly step is ID, i.e., RRQR.  It has less than 5\% of
the flops, but takes 28\% and 37\% time respectively. This shows that
the data movement associated with column pivoting is expensive.

Using the same tolerance $1E-3$, the {\inc} strategy achieves sufficient
accuracy, while the {\dbl} strategy does more adaptations, leading to
higher rank and the acuracy level more than needed, and taking longer time.

\input{tables/tablePerformanceBreakdown.tex}

\subsection{Adaptivity Performance}

\subsubsection{Evaluation of New Stopping Criterion}

This section considers the parametrized problem $\alpha I + \beta UDV^{*}$ with $N=20,000$, $\alpha=1$, $\beta=1$ and $D_{k,k} = 2^{-53(k-1)/r}$, with $r=200$, using $p = 256$ cores.
We provide a comparison with the classical stopping criterion proposed in
HMT \cite{RandomReview2011}, as shown in Table \ref{table:HMT}.
In our strategy, we vary both relative (rtol) and absolute (atol) tolerance.
HMT only works with absolute tolerance, so we put it in the last row of
the table.

It can be seen that with our new criterion, when we set ``rtol'' and ``atol''
to be the same, the approximation accuracy is at the same level of
the requested tolerance (see the diagonal of Table \ref{table:HMT}).

With the HMT criterion, however, since it uses an upper bound as a termination
metric, this in practice can be pessimistic, delivering an approximation error that is usually smaller than requested, at the expense of large ranks. From a user point of view it is difficult to choose a proper compression tolerance.
% This following sentence may or may not be of value.
%This will be further complicated if different subblocks vary greatly
%in magnitude unless one takes care to vary the absolute tolerance
%depending on the magnitude of the elements.

As a result, when the given tolerance is close to machine precision,
HMT criterion requires many steps, as depicted in the bottom right corner of Table \ref{table:HMT}, where the algorithm terminated not due to achieving the ``atol'', but because it reach the maximum allowable rank (5000 in this case).
Thus, there will be some absolute tolerances which
cannot be satisfied, yet it is not clear how small this tolerance is
or how to determine it \emph{before} attempting compression.
Furthermore, it is not clear how to determine when this happens
during compression, either.
A relative tolerance is much easier to set and
frequently of most practical interest.

\input{tables/tableHMT.tex}

\subsubsection{Evaluation of Adaptivity Cost}

In order to assess the cost associated with our adaptivity strategy,
we performed the following experiments with two matrices for which
we know the HSS ranks. 
Thus, we can choose a precise
number of random vectors, so that there is no need for adaptation. This
should be the fastest possible case, and we denote this as ``known-rank''.
Suppose we do not have adaptive strategy, the best a user can do is to
restart compression from scratch (unable to reuse the partial compression
results) when the final residual is large, manually increasing the number of
random vectors. We denote this as ``hard-restart''.
In between these two modes is our adaptive strategy {\inc}.

The first matrix is $I + UDV^*$, where $D_{k,k} = 2^{-53(k-1)/r}, r = 1200, N = 60,000$.
The second matrix is the BEM Acoustic problem, with $N = 10,000$.
Table~\ref{table:adaptiveRank} shows the compression times with different
configurations of $d_0$ and $\Delta d$. It is clear that our adaptive
strategy is nearly as fast as the fastest ``known-rank'' case, and is
up to 2.7x faster than the ``hard-restart''.

\input{tables/tableHardRestart.tex}

\ignore{
  \subsubsection{Restart Overhead}
  Suppose we know the rank before hand, we can do two tests: one using the
  precise number of random vectors without the need for restart, another
  with with different values of initial number of random vectors $d_0$, and
  different $\Delta d$. Or try $d_{i+1} = d_i + \Delta d$ or
  $d_{i+1} = 2 d_i$.
  \textcolor{red}{This can already be seen in the table above.}
}

\subsection{Scalability}

\begin{figure}[H]
	\centering\includegraphics[width=0.4\linewidth]{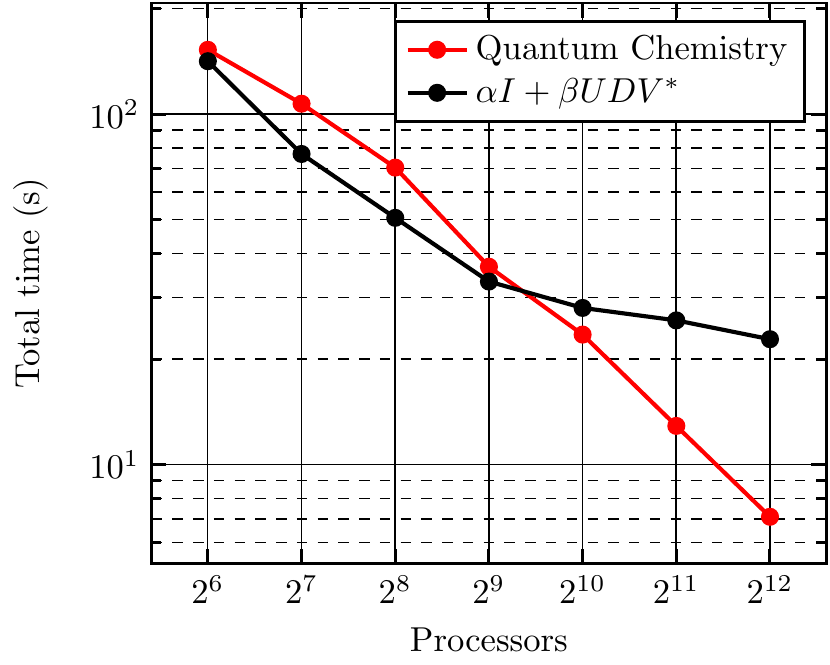}
	\caption{Strong scaling experiment up to $p = 4,096$ cores of the Cori supercomputer. The Quantum Chemistry problem has $N=300,000$ and HSS rank of 12, and the parametrized problem has $N=500,000$ and HSS rank of 480.}
	\label{fig:scaling}
\end{figure}

The last numerical experiment depicts two strong scaling studies, as shown in Figure \ref{fig:scaling}.
The first test problem depicts the Quantum Chemistry dense linear system of size $N=300,000$ and HSS relative approximation error of $1E-2$ and HSS leaf size of 128.
This matrix is amenable to efficient rank compression, resulting in an HSS rank of 12.
In contrast, for a problem with larger numerical rank and tighter numerical accuracy, we show the scalability from the parametrized test case $\alpha I + \beta UDV^{*}$ with $N=500,000$, $\alpha=1$, $\beta=1$ and $D_{k,k} = 2^{-53(k-1)/r}$, with $r=500$.
The resulting HSS rank is $r=480$ at an HSS relative approximation tolerance of $1E-14$ with HSS leaf size of 128. The first example scales better since its HSS rank is small and there is no need for adaptation, whereas the second example requires a few steps for rank adaptation, given that its HSS rank is quite large.

%% file: tables/table6drivers.tex
\begin{table}[H]
\centering%
\caption{Solvings linear systems from different applications.}%
\label{table:sixdrivers}%
\resizebox{\textwidth}{!}{\renewcommand{\arraystretch}{1.1}\setlength{\tabcolsep}{.5em}%
\begin{tabular}{|c|c|c|c|c|c|c|c|c|c|c|c|}
\hline
\multirow{2}{*}{\textbf{Matrix}} & \multirow{2}{*}{$\mathbf{N}$} & \multirow{2}{*}{\textbf{HSS}$_\epsilon$} & \multicolumn{4}{c|}{\textbf{HSS compression}} & \multicolumn{3}{c|}{\textbf{$ULV^*$ factorization}} & \multicolumn{2}{c|}{\textbf{Solve}} \\ \cline{4-12} 
 &  &  & \begin{tabular}[c]{@{}c@{}}HSS\\ rank\end{tabular} & \begin{tabular}[c]{@{}c@{}}Mem\\ (MB)\end{tabular} & \begin{tabular}[c]{@{}c@{}}Flops\\ $\times 10^{12}$\end{tabular} & \begin{tabular}[c]{@{}c@{}}Time\\ (s)\end{tabular} & \begin{tabular}[c]{@{}c@{}}Mem\\ (MB)\end{tabular} & \begin{tabular}[c]{@{}c@{}}Flops\\ $\times 10^{12}$\end{tabular} & \begin{tabular}[c]{@{}c@{}}Time\\ (s)\end{tabular} & \begin{tabular}[c]{@{}c@{}}Flops\\ $\times 10^9$\end{tabular} & \begin{tabular}[c]{@{}c@{}}Time\\ (s)\end{tabular} \\ \hline \hline

\multirow{3}{*}{\textbf{\begin{tabular}[c]{@{}c@{}}$\alpha I+\beta UDV^*$\\(dp)\end{tabular}}} & \multirow{3}{*}{500k}
    & 1E-02 & 28 & 356 & 0.01 & 4.22 & 559  & 0.01 & 0.17 & 0.19 & 0.10  \\ \cline{3-12} 
 &  & 1E-06 & 59 & 671 & 0.03 & 4.71 & 1343 & 0.03 & 0.18 & 0.45 & 0.14  \\ \cline{3-12} 
 &  & 1E-10 & 73 & 868 & 0.06 & 4.93 & 1927 & 0.06 & 0.24 & 0.64 & 0.15  \\ \hline \hline

\multirow{3}{*}{\textbf{\begin{tabular}[c]{@{}c@{}}Toeplitz\\(dp)\end{tabular}}} & \multirow{3}{*}{500k}
    & 1E-02 & 2 & 40 & 0.70 & 1.50 & 64 & 0.001 & 0.11 & 0.02 & 0.06 \\ \cline{3-12} 
 &  & 1E-06 & 2 & 40 & 0.70 & 1.52 & 64 & 0.001 & 0.11 & 0.02 & 0.08 \\ \cline{3-12} 
 &  & 1E-10 & 2 & 40 & 0.70 & 1.62 & 64 & 0.001 & 0.11 & 0.02 & 0.09 \\ \hline \hline

\multirow{3}{*}{\textbf{\begin{tabular}[c]{@{}c@{}}Quantum\\Chem. (dp)\end{tabular}}} & \multirow{3}{*}{500k}
    & 1E-02 & 12 & 235 & 34.60 & 9.26 & 383 & 0.01 & 0.11 & 0.12 & 0.09  \\ \cline{3-12} 
 &  & 1E-06 & 75 & 308 & 34.61 & 9.51 & 486 & 0.01 & 0.22 & 0.16 & 0.09  \\ \cline{3-12} 
 &  & 1E-10 & 113 & 377 & 34.62 & 9.66 & 615 & 0.01 & 0.28 & 0.21 & 0.09 \\ \hline \hline

\multirow{3}{*}{\textbf{\begin{tabular}[c]{@{}c@{}}BEM\\Acoustic (dp)\end{tabular}}} & \multirow{3}{*}{10k}
    & 1E-02 & 723 & 363 & 2.49 & 6.10 & 770 & 0.40 & 0.91 & 0.54 & 0.09     \\\cline{3-12} 
 &  & 1E-06 & 1249 & 607 & 5.12 & 12.72 & 1304 & 1.12 & 1.91 & 0.91 & 0.17  \\ \cline{3-12} 
 &  & 1E-10 & 1332 & 631 & 5.20 & 12.90 & 1366 & 1.23 & 1.97 & 0.95 & 0.27  \\ \hline \hline

 \multirow{3}{*}{\textbf{\begin{tabular}[c]{@{}c@{}}BEM\\Sphere (sp)\end{tabular}}} & \multirow{3}{*}{27k}
    & 1E-01 & 500 & 159 & 9.77 & 11.16 & 364 & 0.13 & 0.55 & 0.47 & 0.09    \\ \cline{3-12} 
 &  & 1E-03 & 1491 & 433 & 26.17 & 29.85 & 1089 & 1.20 & 1.91 & 1.39 & 0.16 \\ \cline{3-12} 
 &  & 1E-05 & 2159 & 795 & 46.97 & 60.07 & 1908 & 3.57 & 4.63 & 2.51 & 0.20 \\ \hline \hline

 \multirow{3}{*}{\textbf{\begin{tabular}[c]{@{}c@{}}Schur100\\(sp)\end{tabular}}} & \multirow{3}{*}{10k}
    & 1E-01 & 287 & 42 & 1.00 & 1.98 & 93 & 0.02 & 0.33 & 0.12 & 0.09  \\ \cline{3-12} 
 &  & 1E-03 & 412 & 81 & 1.20 & 2.60 & 193 & 0.08 & 0.43 & 0.26 & 0.08 \\ \cline{3-12} 
 &  & 1E-05 & 471 & 111 & 1.77 & 3.67 & 257 & 0.14 & 0.44 & 0.34 & 0.04  \\ \hline
\end{tabular}%
}%
\end{table}

%% file: tables/tablePerformanceBreakdown.tex
\begin{table}[h]
\centering
{\footnotesize
\caption{Performance breakdown. ``Compute Samples'' and ``Reduce Samples''
correspond to the two functions in Algorithm~\ref{algo::computelocalsamples}.}
\label{table:PerformanceBreakdown}
\begin{tabular}{l|l|l|l|l|}
\cline{2-5}
 & \multicolumn{2}{c|}{\textbf{Flops $\times 10^{12}$ (\% Flops)}} & \multicolumn{2}{c|}{\textbf{Time (\% Time)}} \\ \cline{2-5} 
 & \inc & \dbl & \inc & \dbl \\ \hline
\multicolumn{1}{|l|}{\textbf{Compression}}        & 29.75 & 56.07 & 5.50 & 8.24 \\ 
\multicolumn{1}{|r|}{HSS Rank}        & \qquad 1480 & \qquad 1945 &  &  \\ 
\multicolumn{1}{|r|}{Relative error}        & \qquad 2.40E-03 & \qquad 2.68E-04 &  &  \\ \hline
\multicolumn{1}{|l|}{$\rightarrow$ Random samp.}  & 26.61 (89.5\%)   & 50.22 (89.6\%) & 2.40 (43.5\%) & 3.93 (47.8\%) \\ \hline
\multicolumn{1}{|l|}{$\rightarrow$ ID}            & 0.85\;\; (2.9\%)     & 2.73\;\; (4.9\%) & 1.54 (27.9\%) & 3.04 (36.9\%) \\ \hline
\multicolumn{1}{|l|}{$\rightarrow$ QR}            & 0.51\;\; (1.7\%)     & \multicolumn{1}{c|}{-} & 0.52 (9.5\%) & \multicolumn{1}{c|}{-} \\ \hline
\multicolumn{1}{|l|}{$\rightarrow$ Orthogonalize} & 0.61\;\; (2.1\%)     & \multicolumn{1}{c|}{-} & 0.18 (3.3\%) & \multicolumn{1}{c|}{-} \\ \hline
\multicolumn{1}{|l|}{$\rightarrow$ Comp. samples} & 1.08\;\; (3.6\%)     & 2.89\;\; (5.2\%) & 0.55 (10.0\%) & 0.86 (10.4\%) \\ \hline
\multicolumn{1}{|l|}{$\rightarrow$ Red. samples}  & 0.09\;\; (0.3\%)     & 0.23\;\; (0.4\%) & 0.32 (5.8\%) & 0.41 (5.0\%) \\ \hline
\multicolumn{1}{|l|}{\textbf{Factorization}}      & 1.29             & 2.67     & 1.98   & 2.83 \\ \hline
\multicolumn{1}{|l|}{\textbf{Solve}}              & 0.001            & 0.002   & 0.14   & 0.15 \\ \hline
\end{tabular}
}
\end{table}

%% file: tables/tableHMT.tex
\begin{table}[]
\centering
{
\footnotesize
\caption{New stopping criterion and HMT criterion \cite{RandomReview2011}.
Each entry has two numbers: one is the relative error of the HSS
approximation given by $||A-HSS*I||_F/||A||_F$, another is the HSS rank.}
\label{table:HMT}
\begin{tabular}{cc|l|l|l|l|}
\cline{3-6}
\multicolumn{1}{l}{}                                 & \multicolumn{1}{l|}{} & \multicolumn{4}{c|}{\textbf{atol}}                                                                                \\ \cline{3-6} 
\multicolumn{1}{l}{}                                 & \multicolumn{1}{l|}{} & \multicolumn{1}{c|}{\textbf{1E-02}} & \multicolumn{1}{c|}{\textbf{1E-06}} & \multicolumn{1}{c|}{\textbf{1E-10}} & \multicolumn{1}{c|}{\textbf{1E-14}} \\ \hline
\multicolumn{1}{|c|}{\multirow{4}{*}{\textbf{rtol}}} & \textbf{1E-02}       &  1.05E-02/43 & 1.05E-02/43  & 1.05E-02/43& 1.05E-02/43 \\ \cline{2-6} 
\multicolumn{1}{|c|}{}                               & \textbf{1E-06}       &  1.82E-05/77 & 1.82E-05/77  & 1.82E-05/77& 1.82E-05/77 \\ \cline{2-6} 
\multicolumn{1}{|c|}{}                               & \textbf{1E-10}       &  4.91E-06/87 & 5.18E-09/127 & 5.18E-09/127& 5.18E-09/127 \\ \cline{2-6} 
\multicolumn{1}{|c|}{}                               & \textbf{1E-14}       &  4.91E-06/87 & 6.68E-10/138 & 6.58E-13/187& 6.58E-13/187 \\ \hline \hline
\multicolumn{1}{l|}{}                                & \textbf{HMT}         &  3.14E-06/87 & 3.50E-10/139 & 3.61E-14/192& 5.53E-15/5000+ \\ \cline{2-6} 
\end{tabular}
}
\end{table}

%% file: tables/tableHardRestart.tex
\begin{table}[]
\centering
{
\footnotesize
\caption{Evaluation of adaptivity cost. For all the tests, we use
 $atol= rtol = 1E-14$, We use $1,024$ cores for the first problem and
64 cores for the second one. With each strategy, we report the
compression time (Compr. time), HSS rank, and the number of adaptation steps
needed (\# adapt.)
 }
\label{table:adaptiveRank}
\begin{tabular}{|l|l|c|c|c|}
\hline
 &  & ``Known-rank'' & {\inc} & ``Hard-restart'' \\ \hline \hline
\multirow{3}{*}{\begin{tabular}[c]{@{}c@{}}$I+UDV^*$\\ $d_0 = 128$\\ $\Delta d = 64$\end{tabular}} & Compr. time & 36.5 & 37.2 & 100.3 \\ \cline{2-5} 
 & HSS-rank & 1162 & 1267 & 1165 \\ \cline{2-5} 
 & Num. adapt. & 0 & 17 & 4 \\ \hline
\multirow{3}{*}{\begin{tabular}[c]{@{}c@{}}$I+UDV^*$ \\ $d_0 = 512$ \\
 $\Delta d = 128$ \end{tabular}} 
 & Compr. time & 35.7 & 36.6 & 54.6\\ \cline{2-5} 
 & HSS-rank & 1162 & 1194 & 1161 \\ \cline{2-5} 
 & \# adapt. & 0 & 5 & 2 \\ \hline \hline
\multirow{3}{*}{\begin{tabular}[c]{@{}c@{}}BEM Acoustics\\ $d_0 = 128$ \\ 
  $\Delta d = 64$ \end{tabular}} 
  & Compr. time & 11.15 & 12.5 & 18.3 \\ \cline{2-5} 
 & HSS-rank & 1264 & 1334 & 1348 \\ \cline{2-5} 
 & \# Adapt. & 0 & 24 & 4 \\ \hline
\multirow{3}{*}{\begin{tabular}[c]{@{}c@{}}BEM Sphere \\ $d_0 = 512$ \\ 
 $\Delta d = 128$ \end{tabular}} 
 & Compr. time & 11.9 & 10.5 & 18.9 \\ \cline{2-5} 
 & HSS-rank & 1276 & 1352 & 1362 \\ \cline{2-5} 
 & \# Adapt. & 0 & 9 & 2 \\ \hline
\end{tabular}
}
\end{table}

%% file: tex/conclusion.tex
%%%%%%%%%%%%%%%%%%%%%%%%%%%%%%%%%%%%%%%%%%%%%%%%%%%%%%%%%%%%%%%%%%%%%%%%
%%% Conclusion
%%%%%%%%%%%%%%%%%%%%%%%%%%%%%%%%%%%%%%%%%%%%%%%%%%%%%%%%%%%%%%%%%%%%%%%%

\section{Conclusion}

We presented two new stopping criteria which allow to accurately
predict the quality of low-rank approximations computed using
randomized sampling. This helps reduce the total number of random
samples as well as reduce the communication cost. We apply these
adaptive randomized sampling schemes for the construction of
hierarchically semi-separable matrices. Compared to previous adaptive
randomized HSS compression approaches, our new methods are more
rigorous and include both absolute and relative stopping criteria. The
numerical examples show faster compression time for the new
incremental adaptive strategy compared to previous methods. Randomized
numerical linear algebra methods are very interesting from a
theoretical standpoint. However, writing robust and efficient,
parallel software is far from trivial. In this paper we have focused
on a number of practical issues regarding the adaptive compression of
HSS matrices, leading to faster compression, with guaranteed
accuracy. The methods shown here should carry over to other structured
matrix representations.

%% file: tex/acknowledgment.tex
\section*{Acknowledgments}
%\todo[inline]{SciDAC, FastMATH, ECP, Yang, Wissam, Liza?}
This research was supported in part by the Exascale Computing
Project (17-SC-20-SC), a collaborative effort of the U.S.
Department of Energy Office of Science
and the National Nuclear Security Administration,
and in part by the U.S. Department of Energy, Office of Science,
Office of Advanced Scientific Computing Research,
Scientific Discovery through Advanced Computing (SciDAC) program.

This research used resources of the National Energy Research Scientific Computing Center (NERSC), a U.S. Department of Energy Office of Science User Facility operated under Contract No. DE-AC02-05CH11231.

We thank Guillaume Sylvand (Airbus) for providing us with the BEM test problems. We thank Daniel Haxton (LBNL) and Jeremiah Jones (Arizona State University) for providing us with the Quantum Chemistry test problem.

We thank Yang Liu (LBNL), Wissam Sid Lakhdar (LBNL), and Liza Rebrova (UCLA) for the insightful
discussions throughout this research.

%%% Local Variables:
%%% mode: latex
%%% TeX-master: t
%%% End:

%% file: tex/appendix.tex
%%%%%%%%%%%%%%%%%%%%%%%%%%%%%%%%%%%%%%%%%%%%%%%%%%%%%%%%%%%%%%%%%%%%%%%%
%%% Appendex for Mathematical Theory

\appendix
\section{Probability Theory}
\label{sec:appendix}
Here, we go through the details of the probability theory mentioned in
Section~\ref{ssec:math_theory}. We define the random variables
\begin{align}
  X &\sim \sigma_{1}^{2}\xi_{1}^{2} + \cdots + \sigma_{r}^{2}\xi_{r}^{2}
      \nonumber\\
  \overline{X}_{d} &\sim \frac{1}{d}\brackets{X_{1} + \cdots + X_{d}},
\end{align}
where $\xi_{j} \sim \mathcal{N}(0,1)$,
$\sigma_{1}\ge\cdots\ge\sigma_{r}>0$ are the singular values of $A$,
and $X_{i}$ are independent realizations of $X$.  Clearly
$\E\parens{X} = \E\parens{\overline{X}_{d}} = \norm{A}_{F}^{2}$.

We prove the following probabilistic bounds:
\begin{align}
  \label{eq:A_prob_bounds}
  \P\brackets{\overline{X}_{d}\ge \norm{A}_{F}^{2}\tau}
  &\le \exp\parens{-\frac{d\tau}{2}} \norm{A}_{F}^{dr}
    \prod_{k=1}^{r}\parens{A_{k}'}^{-d}
    \qquad \tau>1 \nonumber\\
  \P\brackets{\overline{X}_{d}\le \norm{A}_{F}^{2}\tau}
        &\le \exp\parens{\frac{d\tau}{2}} \norm{A}_{F}^{dr}
          \prod_{k=1}^{r}\parens{A_{k}''}^{-d}
          \qquad \tau\in[0,1) \, ,
\end{align}
where
\begin{align}
  \norm{A}_{F}^{2} &= \sigma_{1}^{2} + \cdots + \sigma_{r}^{2} \nonumber\\
  \parens{A_{k}'}^{2} &= \norm{A}_{F}^{2} - \sigma_{k}^{2}  \label{eq:As_def} \\
    \parens{A_{k}''}^{2} &= \norm{A}_{F}^{2} + \sigma_{k}^{2} \, . \nonumber
\end{align}
We will show particular realizations of the random variable deviate
from the mean with exponentially-decaying probability and gives an
accurate way to estimate $\norm{A}_{F}$ using random (Gaussian)
samples of the range. Thus, the new stopping criterion accurately
estimates the error; the singular values determine the exact
exponential-decaying probability.

\subsection{Chernoff's Inequality}
We use the following theorem~\cite{prob_book}:
\begin{theorem}[Chernoff's Inequality]\label{thm:Chernoff}
  Given a random variable $X$, we have
  \begin{equation}
    \label{eq:Chernoff}
    \P\brackets{X \ge a} \le \min_{t > 0} e^{-ta}\, \E\parens{e^{tX}} \, .
  \end{equation}
\end{theorem}
Here, $\E\parens{e^{tX}}$ is the moment generating function
of a random variable $X$.
A slight modification of Theorem~\ref{thm:Chernoff}
gives
\begin{equation}
  \label{eq:Chernoff_2}
  \P\brackets{X \le a} \le \min_{t>0} e^{ta}\, \E\parens{e^{-tX}} \, .
\end{equation}
$X$ is merely linear combinations of chi-squared distributions,
so by using properties of the moment generating function we see
%We see \todo{where does this come from? I do not see}
\begin{equation}
  M_{\overline{X}_{d}}(t) =
  \prod_{k=1}^{d}
  \parens{1 - \frac{2\sigma_{k}^{2}}{d}t}^{-\frac{d}{2}}.
\end{equation}
Unfortunately, it is difficult to compute
\begin{equation}
  \min_{t>0} e^{-ta}M_{\overline{X}_{d}}(t)
\end{equation}
analytically, as it reduces to computing the zeros of a degree $r$
polynomial. Instead, we settle for choosing a particular value
of $t$. Setting
\begin{equation}
  \bar{t} = \frac{d}{2\norm{A}_{F}^{2}}
\end{equation}
and $a = \norm{A}_{F}^{2}\tau$, with $\tau>1$, we find
\begin{equation}
  \P\brackets{\overline{X}_{d}\ge \norm{A}_{F}^{2}\tau}
  \le \exp\parens{-\frac{d\tau}{2}} \norm{A}_{F}^{dr}
  \prod_{k=1}^{r}\parens{A_{k}'}^{-d}
  \qquad \tau>1,
\end{equation}
where $A_{k}'$ is defined in Eq.~\eqref{eq:As_def}. In a similar
manner, we can obtain a lower bound
\begin{equation}
  \P\brackets{\overline{X}_{d}\le \norm{A}_{F}^{2}\tau}
  \le \exp\parens{\frac{d\tau}{2}} \norm{A}_{F}^{dr}
  \prod_{k=1}^{r}\parens{A_{k}''}^{-d}
  \qquad \tau\in[0,1) \, .
\end{equation}
It is clear that these lower bounds are not optimal and more rigorous
analysis would produce tighter bounds, but we will not investigate
this further. The bounds we just obtained are sufficient for our
purposes.

\subsection{Exponential Decay}
We now work through the details to show we have exponential-decaying
tail probabilities once we are far enough away from the expected
value.  The restrictions are mild: for lower bounds, we require
$\tau < \ln2$; for upper bounds we require
$\tau > 1 + \frac{\norm{A}_{2}^{2}}{\norm{A}_{F}^{2} -
  \norm{A}_{2}^{2}}$. The probability distributions
do not apply for rank-1 matrices (when $\norm{A}_{2} = \norm{A}_{F}$),
but in then Chernoff's Inequality can be used to
compute optimal bounds on tail probabilities. Therefore,
we ignore this case.

%The only time we have
%$\norm{A}_{F} = \norm{A}_{2}$ is when $A$ is a rank-1 matrix, in which
%case we can compute the optimal bounds using Chernoff's inequality.
%Furthermore, because $\E\parens{X} = \norm{A}_{F}^{2}$, we will take a
%square root to approximate $\norm{A}_{F}$, further dampening
%deviations.  As long as we are concerned with having a value that is
%the correct order of magnitude, this method should be
%sufficient.
%\todo{remove this comment?}; Removed -- Chris

\subsubsection{$\P\brackets{\overline{X}_{d}\ge \norm{A}_{F}^{2}\tau}$}

For $\tau>1$, we have
\begin{align}
    \P\brackets{\overline{X}_{d}\ge \norm{A}_{F}^{2}\tau}
        &\le \exp\parens{-\frac{d\tau}{2}} \norm{A}_{F}^{dr}
            \prod_{k=1}^{r}\parens{A_{k}'}^{-d} \nonumber\\
        &= \exp\brackets{\frac{d}{2}\braces{\parens{\nu_{1}
            + \cdots + \nu_{r}} - \tau}} \, ,
\end{align}
where
\begin{equation}
    \nu_{k} = \ln\brackets{\frac{1}{
        1 - \frac{\sigma_{k}^{2}}{\norm{A}_{F}^{2}}}} \, .
\end{equation}
To have exponential decay in probability, we require
\begin{equation}
    \nu_{1} + \cdots + \nu_{r} < \tau \, .
\end{equation}
Because $-\ln x$ is convex on $\parens{0,\infty}$, we see that for
$\alpha\in\parens{0,1}$, we have
\begin{equation}
    \ln\parens{\frac{1}{1-x}} \le \frac{x}{\alpha}\ln\parens{\frac{1}{1-\alpha}}
        \qquad x\in\brackets{0,\alpha} \, .
\end{equation}
Now, we know
$\frac{\sigma_{k}^{2}}{\norm{A}_{F}^{2}} \le
\frac{\norm{A}_{2}^{2}}{\norm{A}_{F}^{2}}$, so we let
\begin{equation}
  \label{eq:alpha_val}
  \alpha = \frac{\norm{A}_{2}^{2}}{\norm{A}_{F}^{2}}.
\end{equation}
Combing this with the fact that $\ln(1+x)\le x$, we see
\begin{align}
    \nu_{1} + \cdots + \nu_{r}
        &\le \frac{\sigma_{1}^{2}}{\norm{A}_{F}^{2}}\frac{1}{\alpha}
            \ln\parens{\frac{1}{1-\alpha}} + \cdots +
            \frac{\sigma_{r}^{2}}{\norm{A}_{F}^{2}}\frac{1}{\alpha}
            \ln\parens{\frac{1}{1-\alpha}} \nonumber\\
        &\le \frac{1}{1-\alpha} \, .
\end{align}
If we plug in $\alpha$ from Eq.~\eqref{eq:alpha_val}, we see that we
require
\begin{equation}
    \tau > 1 + \frac{\norm{A}_{2}^{2}}{\norm{A}_{F}^{2} - \norm{A}_{2}^{2}}
\end{equation}
in order to have exponentially decaying tail probabilities in $d$.

\subsubsection{$\P\brackets{\overline{X}_{d}\le \norm{A}_{F}^{2}\tau}$}

For $\tau\in[0,1)$, we have
\begin{align}
    \P\brackets{\overline{X}_{d}\le \norm{A}_{F}^{2}\tau}
        &\le \exp\parens{\frac{d\tau}{2}} \norm{A}_{F}^{dr}
            \prod_{k=1}^{r}\parens{A_{k}''}^{-d} \nonumber\\
        &= \exp\brackets{\frac{d}{2}\braces{\tau - \parens{\lambda_{1}
            + \cdots + \lambda_{r}}}} \, ,
\end{align}
where
\begin{equation}
    \lambda_{k} = \ln\brackets{1 + \frac{\sigma_{k}^{2}}{\norm{A}_{F}^{2}}}.
\end{equation}
To have exponential decay in probability, we require
\begin{equation}
    \lambda_{1} + \cdots + \lambda_{r} > \tau \, .
\end{equation}
Now, we know
\begin{equation}
    \ln\parens{1 + x} \ge x\ln 2 \qquad x\in\brackets{0,1},
\end{equation}
so this implies
\begin{align}
    \lambda_{1} + \cdots + \lambda_{r}
        &\ge \frac{\sigma_{1}^{2}}{\norm{A}_{F}^{2}}\ln2
            + \cdots + \frac{\sigma_{r}^{2}}{\norm{A}_{F}^{2}}\ln2 \nonumber\\
        &= \ln 2.
\end{align}
Therefore, so long as $\tau<\ln2$, we have exponentially decaying tail
probabilities in $d$.